\documentclass[11pt]{article}
\usepackage[table]{xcolor}
\usepackage{graphicx}
\usepackage{amsmath,amsfonts,amssymb}
\usepackage{graphicx}
\usepackage{multirow}
\usepackage{tikz,pgf}
\usepackage{url}
\usepackage{amsmath}
\usepackage{graphicx}
\usepackage{epsfig}
\usepackage{epstopdf}
\usepackage{color}
\usepackage{enumitem}
\def\tto{\;{\lower 1pt \hbox{$\rightarrow$}}\kern -10pt
\hbox{\raise 2pt \hbox{$\rightarrow$}}\;}

\def\Tilde{\widetilde}
\def\Bar{\overline}
\def\ra{\rangle}
\def\la{\langle}

\def\epsilon{\varepsilon}
\def\B{\Bbb B}
\def\h{\hfill\Box}
\def\R{\Bbb R}

\def\N{\Bbb N}
\def\ox{\bar{x}}

\def\oy{\bar{y}}

\def\co{\mbox{\rm co}\,}

\def\ri{\mbox{\rm ri}\,}

\def\gph{\mbox{\rm gph}\,}
\def\aff{\mbox{\rm aff}\,}
\def\epi{\mbox{\rm epi}\,}

\def\dom{\mbox{\rm dom}\,}
\def\aff{\mbox{\rm aff}\,}

\def\rge{\mbox{\rm rge}\,}

\def\h{\hfill\square}

\def\ph{\varphi}

\def\oR{\Bar{\R}}
\def\s{\square}

\def\ph{\varphi}

\def\oR{\Bar{\R}}

\setlist[enumerate,1]{itemsep=0.0ex,parsep=0.5ex,label={\rm(\alph*)},leftmargin=*, align=left}
\newcounter{lk}

\usepackage[colorlinks=true]{hyperref}

\begin{document}
\begin{center}
{\sc\bf  Nearly Convex Optimal Value Functions\\   and Some Related Topics}\\[1ex]
{{\sc Nguyen Quang Huy}\footnote{Department of Science Management and International Relations, Hanoi Pedagogical University 2, Vinh Phuc, Vietnam; email: nqhuy@hpu2.edu.vn; huyngq308@gmail.com.},
{\sc Nguyen Mau   Nam}\footnote{Fariborz Maseeh Department of Mathematics and Statistics, Portland State University, Portland, OR
97207, USA (mnn3@pdx.edu). Research of this author was partly supported by the USA National Science Foundation under grant DMS-2136228.},
{\sc Nguyen Dong Yen}\footnote{Institute of Mathematics, Vietnam Academy of Science and Technology, 18 Hoang Quoc Viet,
Hanoi 10307 (ndyen@math.ac.vn).}}\\
\end{center}
\small{\bf Abstract.} In this paper, we introduce new properties of the relative interior calculus for nearly convex sets, functions, and set-valued mappings. These properties are important for the development of duality theory in optimization. Then we  investigate optimal value functions defined by nearly convex functions and nearly convex set-valued mappings, and derive the near convexity of the optimal value function under a qualification condition. We also develop formulas for subgradients and Fenchel conjugates of this class of functions, and explore their applications to duality theory. \\[1ex]

\newtheorem{Theorem}{Theorem}[section]
\newtheorem{Proposition}[Theorem]{Proposition}
\newtheorem{Remark}[Theorem]{Remark}
\newtheorem{Lemma}[Theorem]{Lemma}
\newtheorem{Corollary}[Theorem]{Corollary}
\newtheorem{Definition}[Theorem]{Definition}
\newtheorem{Example}[Theorem]{Example}
\renewcommand{\theequation}{\thesection.\arabic{equation}}
\normalsize

\section{Introduction}
\setcounter{equation}{0}

Convex analysis has been extensively studied since the late 20th century, with significant contributions from many mathematicians. Its important role in convex optimization, along with recent applications to fields such as economics, computational statistics, and machine learning, has motivated the search for new notions of generalized convexity. In the Euclidean space $\R^n$, Minty and Rockafellar introduced the concept of {\em near convexity} for sets in the 1960s, which requires that the set under consideration lies between a convex set and its closure (see~\cite{Minty1961, R1970}). As a natural generalization of convexity, there are numerous examples of nearly convex sets that are not convex. For instance, removing any non-vertex point on a side of a square in two dimensions results in a nearly convex set that is not convex. The near convexity is also interesting because it appears in the study of maximally monotone operators due to the fact that the domain of a maximally monotone operator is nearly convex. In particular, the domain of the subdifferential mapping of a proper lower semicontinuous convex function is nearly convex.

The concept of near convexity provides a framework for defining nearly convex functions and set-valued mappings, which have found numerous applications in optimization and duality theory. Recent research, including \cite{bmw2013,bkw2008, LM2019,mmw2016, NTY2023}, has explored the properties of nearly convex sets and functions, as well as their connections to ranges of maximally monotone operators. Among these papers, we  mention the work of Bo{\c{t}}, Kassay, and Wanka from \cite{bkw2008} on properties of nearly convex sets and functions with application to duality theory. In \cite{bmw2013, mmw2016}, further properties of nearly convex sets and their connections to maximally monotone operators were explored. Our recent paper \cite{NTY2023} introduced the notion of nearly convex set-valued mappings and provided a unified study of nearly convex sets, functions, and set-valued mappings. There have been further studies and applications of near convexity in recent literature, including \cite{nr1, nr2, nr3, nr4, nr5}.

In this paper, we continue our work on nearly convex set-valued mappings, aiming at applications to duality theory for more general optimization models than those considered in \cite{bkw2008}. We studied near convexity and subdifferentiation of the optimal value function, which has been known to be very important in the study of parametric optimization. Another notable contribution is a new formula for the Fenchel conjugate of the optimal value function defined by nearly convex objective functions and nearly convex constraint set-valued mappings. This formula is obtain using the Fenchel conjugate of the objective function and the {\em Fenchel conjugate} of the set-valued mapping which we introduce for the first time. Along with the new results obtained, we use a geometric approach to significantly simplify most of the proofs for the results in \cite{bkw2008} on nearly convex sets, functions, and duality.

Our paper is organized as follows. In Section 2, we provide an overview of the fundamental definitions and properties of nearly convex sets, functions, and set-valued mappings. In Section 3 and 4, we present further studies on nearly convex sets, functions, and set-valued mappings that are essential for our later development on duality theory. In Section 5, we establish the near convexity of the optimal value function and derive a new formula for subgradients of this function with near convexity data. Section 6 introduces the concept of Fenchel conjugate of set-valued mappings, where we develop a new formula for the Fenchel conjugate of the optimal value function and its applications to Fenchel conjugate calculus for nearly convex functions. Finally, in Section 7, we explore the applications of our findings to duality theory.

The paper utilizes standard notions and notations of convex analysis in the Euclidean space $\R^n$, which can be found in sources such as \cite{Borwein2000,HU2,bmn,r}. In the sequel, we use the notation $B(z,\rho)$ (resp., $\B (z,\rho)$) to represent the open (resp., closed) ball centered at $z\in\mathbb R^k$ with a radius of $\rho>0$. The affine hull of a subset $D\subset\mathbb R^k$ is abbreviated to $\aff D$, while $\overline{D}$ represents the closure of $D$, and $\mbox{\rm int}\,D$ denotes the interior of $D$.

\section{Preliminaries on Near Convexity}
\setcounter{equation}{0}

This section provides a brief overview of fundamental definitions and significant properties of nearly convex sets that are utilized in this paper. For more comprehensive information, we recommend referring to~\cite{bmw2013,bkw2008, LM2019,mmw2016, NTY2023}.

Let $\Omega$ be a subset of $\R^n$. We say that $\Omega$ is {\em nearly convex} if there exists a convex set $C$ such that
\begin{equation}\label{NC}
C\subset\Omega\subset \Bar{C}.
\end{equation}

Recall that the {\em relative interior} of an arbitrary set $\Omega$  in $\R^n$ is defined by
\begin{equation*}
\ri\Omega=\big\{a\in \Omega\; \big|\; \mbox{\rm there exists }\delta>0\; \mbox{\rm such that }B(a; \delta)\cap \aff\Omega\subset \Omega\big\}.
\end{equation*}
It follows from the definition that $a\in \ri\Omega$ if and only if $a\in \aff\Omega$ and there exists $\delta>0$ such that
\begin{equation*}
B(a; \delta)\cap \aff\Omega\subset \Omega.
\end{equation*}
Since $\Bar{\ri C}=\Bar{C}$ for a convex set $C$, it follows from the definition that $\Omega$ is nearly convex if and only if there exists a convex set $C$ such that
\begin{equation*}
\ri C\subset \Omega\subset \Bar{C}.
\end{equation*}
The proposition below shows that the notion of near convexity coincides with the notion {\em almost convexity} considered in~\cite{bkw2008}.
\begin{Proposition} Let $\Omega$ be a subset of $\R^n$. Then $\Omega$ is nearly convex if and only if $\Bar{\Omega}$ is convex and $\ri \Bar{\Omega}\subset\Omega$.
\end{Proposition}
{\bf Proof.} Suppose that $\Omega$ is nearly convex. Then there exists a convex set $C$ such that~\eqref{NC} is satisfied. It is not hard to see that $\Bar{\Omega}=\Bar{C}$ and $\ri \Bar\Omega=\ri C$. Since $C$ is convex, the set $\Bar{\Omega}=\Bar{C}$ is  convex. In addition, $\ri \Bar\Omega=\ri C\subset C\subset \Omega$. This justifies the implication $\Longrightarrow$ of the proposition.

To prove the reverse implication, suppose that $\Bar{\Omega}$ is convex and $\ri \Bar{\Omega}\subset\Omega$. Let $C=\ri \Bar{\Omega}$. Since $\Bar{\Omega}$ is convex, we see that $C$ is convex. By the given assumption, $C=\ri \Bar{\Omega}\subset\Omega$. We also have by the convexity of $\Bar{\Omega}$ that
\begin{equation*}
\Bar{C}=\Bar{\ri \Bar{\Omega}}=\mbox{\rm cl}\, \Bar{\Omega}=\Bar{\Omega}\supset\Omega.
\end{equation*}
This justifies the near convexity of $\Omega$. $\h$

The theorem below lists some important properties of nearly convex sets; see, e.g., \cite{bkw2008,mmw2016}.
\begin{Theorem}\label{T} The following assertions hold:
\begin{enumerate}
\item If $\Omega_1$ is a nearly convex set in $\R^n$ and $\Omega_2$ is a nearly convex set in $\R^p$, then $\Omega_1\times \Omega_2$ is a nearly convex set in $\R^{n}\times \R^p$.
    \item Let $\Omega$ be a nearly convex set with $C\subset \Omega\subset \Bar{C}$, where $C$ is a convex set in $\R^n$. Then $\aff\Omega=\aff C$, $\ri C=\ri\Omega$, and $\overline{C}=\overline{\Omega}$.
\item Let $\Omega_1,\dots, \Omega_m$ be nearly convex sets in $\R^n$. If $\bigcap\limits_{i=1}^m\ri\Omega_i\neq\emptyset$, then $\bigcap\limits_{i=1}^m\Omega_i$ is nearly convex and
\begin{equation*}
\ri\left(\bigcap_{i=1}^m\Omega_i\right)=\bigcap_{i=1}^m\ri\Omega_i.
\end{equation*}
\item If $\Omega$ is a nearly convex set in $\R^n$ and $A\colon \R^n\to \R^p$ is a linear mapping, then $A(\Omega)$ is a nearly convex set in $\R^p$ and $\ri A(\Omega)=A(\ri\Omega)$.
\end{enumerate}
\end{Theorem}

Given a function $f\colon \R^n\to \oR=[-\infty, \infty]$, the {\em effective domain} and the  {\em epigraph} of $f$ are given  respectively by
\begin{align*}
&\dom f=\{x\in \R^n\; |\; f(x)<\infty\},\\
&\epi f=\{(x, \lambda)\in \R^n\times \R\; |\; f(x)\leq \lambda\}.
\end{align*}
The function $f$ is said to be {\em convex} if $\epi f$ is a convex set, and it is said to be {\em nearly convex} if $\epi f$ is a nearly convex set. We say that $f$ is {\em proper} if $\dom f\neq\emptyset$ and $-\infty <f(x)$ for all $x\in \R^n$.

For a set-valued mapping $F\colon \R^n\tto \R^p$, define the {\em domain},  the \textit{range}, and the {\em graph} of $F$ by
\begin{align*}
&\dom F=\{x\in \R^n\; |\; F(x)\neq\emptyset\},\quad  \rge F=\bigcup\limits_{x\in\R^n}F(x),\\
&\gph F=\{(x, y)\in \R^n\times \R^p\; |\; y\in F(x)\}.
\end{align*}
The set-valued mapping $F$ is said to be {\em proper} if $\dom F\neq\emptyset$. We say that $F$ is  {\em nearly convex} if $\gph F$ is a nearly convex set in $\R^n\times \R^p$.

Given a  function $f\colon \R^n\to \oR$, define the {\em epigraphical mapping} $E_f\colon \R^n\tto \R$ by
\begin{equation}\label{Emapping}
   E_f(x)=\big\{\lambda\in \R\;\big|\; f(x)\leq \lambda\big\},\ \; x\in \R^n.
\end{equation}
It follows directly from the definition that $\dom E_f=\dom f$ and $\gph E_f=\epi f$.

Let us now present an important theorem on the relative interior of the graph of a nearly convex set-valued mapping~\cite[Theorems~3.6 and~3.9]{NTY2023}.

\begin{Theorem}\label{Rthm} Let $F\colon \R^n\tto \R^p$ be a nearly convex set-valued mapping. Then
\begin{equation*}
\mbox{\rm ri}(\gph F)=\big\{(x, y)\in \R^n\times \R^p\; \big|\; x\in \mbox{\rm ri}(\dom F), \; y\in \ri F(x)\big\}.
\end{equation*}
In addition, $F(x)$ is nearly convex  and ${\rm ri}\,F(x)$ is nonempty for every $x\in {\rm ri}(\dom F)$.
\end{Theorem}

Let $g\colon \R^n\to \R^p$ and let $M$ be a nonempty subset of $\R^p$. Following Bo{\c{t} et al.~\cite{bkw2008}, we define the {\em epigraph} of~$g$ with respect to~$M$ by
\begin{equation*}
\mbox{\rm epi}_M(g)=\{(x, y)\in \R^n\times \R^p\; |\; y-g(x)\in M\},
\end{equation*}
We say that $g$ is {\em $M$-nearly convex} if $\mbox{\rm epi}_M(g)$ is nearly convex.  Thanks to Theorem~\ref{Rthm}, we can give an alternative simple proof for the next significant result in~\cite{bkw2008}.

\begin{Theorem}\label{T1}{\rm\bf (\cite[Lemma~2.5]{bkw2008})} Let $g\colon \R^n\to \R^p$ be  $M$-nearly convex. Then
\begin{equation}\label{rep1}
\mbox{\rm ri}({\rm epi}_M(g))=\big\{(x, y)\in \R^n\times \R^p\; |\; y-g(x)\in \ri M\big\}.
\end{equation}
\end{Theorem}
{\bf Proof.} Define the set-valued mapping $F\colon \R^n\tto \R^p$ by $F(x)=g(x)+M$. Since $g$ is $M$-nearly convex,  $\gph F=\big\{(x,y)\in\R^n\times\R^p\mid y\in F(x)\big\}={\rm epi}_M(g)$ is nearly convex by definition. Obviously, $\dom F=\R^n$ and for any $x\in \R^n$ we have $\ri F(x)=g(x)+\ri M$. So, by Theorem~\ref{Rthm}, we have
\begin{equation*}
\mbox{\rm ri}({\rm epi}_M(g))=\mbox{\rm ri}(\gph F)=\big\{(x, y)\in \R^n\times \R^p\; |\; x\in \mbox{\rm ri}(\dom F), \; y \in \ri F(x)\big\}.
\end{equation*}
This immediately implies~\eqref{rep1}. $\h$

Throughout the paper we deal mostly with proper functions but occasionally encounter improper ones. The proposition below shows that improper nearly convex functions are rare.

\begin{Proposition}\label{PP} Let $f\colon \R^n\to\oR$ be a  nearly convex function. If there exists a point $\bar x\in {\rm ri}(\dom f)$ with $f(\bar x)\in\R$, then $f$ is proper, i.e., $f(x)>-\infty$ for all $x\in\R^n$.
\end{Proposition}
{\bf Proof.} On the contrary, suppose that there is a point $\hat x$ such that $f(\hat x)=-\infty$. Then, $(\hat x,\gamma)\in {\rm epi} f$ for every $\gamma\in\R$. Since ${\rm epi} f$ is nearly convex, one has $C\subset {\rm epi} f\subset\overline{C}$ for some convex set $C\subset \R^n\times\R$. Setting ${\mathcal P}_1(x,\alpha)=x$ for $(x,\alpha)\in \R^n\times\R$, one sees that $\dom f ={\mathcal P}_1({\rm epi} f)$. So, by Theorem~\ref{T}(c), we deduce from the near convexity of ${\rm epi} f$ that $\dom f$ is nearly convex. Since $\bar x\in {\rm ri}(\dom f)$ and $\hat x\in\dom f$, taking $\tau>0$ small enough, we have that the point $x_\tau:=\ox +\tau(\ox -\hat x)$ belongs to $\dom f$. Take any $\beta>f(x_\tau)$ and note that $(x_\tau,\beta)\in {\rm epi} f$. Then, by the inclusion ${\rm epi} f\subset\overline{C}$ and the convexity of $\overline{C}$, the line segment joining the points $(x_\tau,\beta)$ and $(\hat x,\gamma)$, where $\gamma\in\R$ can be chosen arbitrarily, is contained in $\overline{C}$. In particular, the point
$$\dfrac{1}{1+\tau}(x_\tau,\beta)+\dfrac{\tau}{1+\tau}(\hat x,\gamma)=\left(\ox,\dfrac{1}{1+\tau}\beta+\dfrac{\tau}{1+\tau}\gamma\right)$$ belongs to $\overline{C}$. Applying Theorem~\ref{Rthm} to $F:=E_f$, where $E_f$ is the epigraphical mapping defined in~\eqref{Emapping} yields $(\bar x,f(\bar x))\notin {\rm ri}({\rm epi} f)$ and $(\bar x,\alpha)\in {\rm ri}({\rm epi} f)$ for any $\alpha>f(\bar x)$. Then, from the equality ${\rm ri}({\rm epi} f)={\rm ri}(\overline{C})$, which is valid by Theorem~\ref{T}(b), and the inclusion $\left(\ox,\dfrac{1}{1+\tau}\beta+\dfrac{\tau}{1+\tau}\gamma\right)\in  \overline{C}$ we can deduce that
$$\left[(\bar x,\alpha),\left(\ox,\dfrac{1}{1+\tau}\beta+\dfrac{\tau}{1+\tau}\gamma\right)\right)\in {\rm ri}(\overline{C}).$$ Noting that the point $(\bar x,f(\bar x))$ is contained in the half-open line segment on the right-hand side of the last inclusion if $\gamma\in\R$ is chosen small enough, we get $(\bar x,f(\bar x))\in {\rm ri}(\overline{C})={\rm ri}({\rm epi} f)$. We have arrived at a contradiction, which completes the proof. $\h$

\section{Relative Interior Calculus under Near Convexity}
\setcounter{equation}{0}

In this section, we study the relative interiors of nearly convex sets and nearly convex set-valued mappings under basic operations on set-valued mappings.

Let us establish the first new result on the relative interior of the direct image of a set under a  nearly convex set-valued mapping.
\begin{Theorem}\label{ri_of_image} Let $F\colon \R^n\tto \R^p$ be a nearly convex set-valued mapping and let $\Omega$ be a nearly convex set in $\R^n$. Suppose that
\begin{equation}\label{Reg1}
\mbox{\rm ri}(\dom F)\cap \ri \Omega\neq\emptyset.
\end{equation}
Then
\begin{equation}\label{R1}
\mbox{\rm ri}\big(F(\Omega)\big)=\bigcup_{x\in (\ri\Omega)\cap \mbox{\rm ri}(\dom F)}\ri F(x).
\end{equation}
\end{Theorem}
{\bf Proof.} Define the nearly convex sets
\begin{equation*}
\Omega_1=\gph F\ \; \mbox{\rm and }\; \Omega_2=\Omega\times \R^p.
\end{equation*}
By~\eqref{Reg1}, there exists a vector $\hat{x}\in \mbox{\rm ri}(\dom F)\cap \ri \Omega$. Since $F(\hat{x})$ is nearly convex by~\cite[Theorem 3.9]{NTY2023}, we can find an element $\hat{y}\in \ri F(\hat{x})$. Using Theorem~\ref{Rthm}, we see that $$(\hat{x}, \hat{y})\in \mbox{\rm ri}(\gph F)=\ri\Omega_1.$$ Since $(\hat{x}, \hat{y})\in\ri\Omega_2$, this yields $\ri\Omega_1\cap \ri\Omega_2\neq\emptyset$; so $\Omega_1\cap \Omega_2$ is nearly convex by~\cite[Theorem~2.1]{bkw2008} (see also~\cite[Corollary~4.8]{mmw2016})). Observe that $F(\Omega)=\mathcal{P}_2(\Omega_1\cap \Omega_2)$, where $\mathcal{P}_2(x,y)=y$ for every $(x,y)\in\R^n\times\R^p$. Then by the linearity of the operator  $\mathcal{P}_2$ and by~\cite[Theorem~4.2 and Corollary~4.8]{mmw2016} one has
\begin{equation}\label{riF_Omega}
\mbox{\rm ri}\big(F(\Omega)\big)=\mathcal{P}_2(\mbox{\rm ri}(\Omega_1\cap \Omega_2))=\mathcal{P}_2(\ri\Omega_1\cap \ri\Omega_2).
\end{equation}
Note that $\ri\Omega_2=(\ri\Omega)\times \R^p$. In addition, by Theorem~\ref{Rthm} we have
\begin{equation*}
\ri \Omega_1=\big\{(x, y)\in \R^n\times \R^p\; |\; x\in \mbox{\rm ri}(\dom F), \; y\in \ri F(x)\big\}.
\end{equation*}
These representations of $\ri\Omega_1$ and $\ri\Omega_2$ together with~\eqref{riF_Omega} clearly imply~\eqref{R1}. $\h$

\medskip
As a corollary of Theorem~\ref{ri_of_image}, we derive a result given by Bo{\c{t}}, Kassay  and Wanka in~\cite{bkw2008}.

\begin{Corollary}\label{ri epi f}{\rm\bf (\cite[Lemma~2.6]{bkw2008})} Suppose that $X\subset \R^n$ is nearly convex and $g\colon \R^n\to \R^p$ is $M$-nearly convex. Then
\begin{equation*}
\mbox{\rm ri}(g(X)+M)=g(\ri X)+\ri M.
\end{equation*}
\end{Corollary}
{\bf Proof.} Consider the set-valued mapping $F$ defined in  the proof of Theorem \ref{T1}. Then $\dom F=\R^n$, and $\ri F(x)=g(x)+\ri M$ for any $x\in \R^n$. It follows from~\eqref{R1} with $\Omega=X$ that
\begin{equation*}
\mbox{\rm ri}\big(F(\Omega)\big)=\bigcup_{x\in \ri X}\ri F(x)=\bigcup_{x\in\ri X}(g(x)+\ri M=g(\ri X)+\ri M.
\end{equation*}
This completes the proof. $\h$

Given a set-valued mapping $F\colon \R^n\tto \R^p$ and a set $\Theta\subset\R^p$, define  the {\em inverse image} of $\Theta$ under $F$ by
\begin{equation*}
F^{-1}(\Theta)=\big\{x\in \R^n\; |\; F(x)\cap \Theta\neq\emptyset\big\}.
\end{equation*}
The next theorem allows us to get a representation of the relative interior of the inverse image of a set under a  nearly convex set-valued mapping.
\begin{Theorem} \label{RINV}
Let $F\colon \R^n\tto \R^p$ be a nearly convex set-valued mapping and let $\Theta$ be a nearly convex set in  $\R^p$. Suppose that
\begin{equation}\label{QC3}
\mbox{\rm ri}(\rge F)\cap \ri \Theta\neq\emptyset.
\end{equation}
Then define $F_0(x)=\ri F(x)$ for $x\in \R^n$ and get
\begin{equation}\label{R2}
\mbox{\rm ri}\big(F^{-1}(\Theta)\big)=\mbox{\rm ri}(\dom F)\cap F_0^{-1}(\ri \Theta).
\end{equation}
\end{Theorem}
{\bf Proof.} Consider the nearly convex sets
\begin{equation*}
\Omega_1=\gph F\ \; \mbox{\rm and }\; \Omega_2=\R^n\times \Theta.
\end{equation*}
 By~\eqref{QC3}, we can find a vector $\hat{y}\in \mbox{\rm ri}(\rge F)\cap \ri \Theta$. Considering  the set-valued mapping $G=F^{-1}$, which is nearly convex, we see that $\hat{y}\in \mbox{\rm ri}(\rge F)=\mbox{\rm ri}(\dom G)$.  So, according to Theorem~\ref{Rthm}, there is a point $\hat{x}\in \ri G(\hat{y})$. Now consider the linear function
\begin{equation*}
T(y, x)=(x, y), \ \; (y, x)\in \R^p\times \R^n.
\end{equation*}
Then we have $T(\gph G)=\gph F$ and thus $T(\mbox{\rm ri}(\gph G))=\mbox{\rm ri}(\gph F)$. Since  the inclusion $(\hat{y}, \hat{x})\in \mbox{\rm ri}(\gph G)$ holds by Theorem~\ref{Rthm}, we see that $(\hat{x}, \hat{y})\in \mbox{\rm ri}(\gph F)=\ri\Omega_1$. Thus, $(\hat{x}, \hat{y})\in \ri \Omega_1\cap \ri\Omega_2$, so $\Omega_1\cap \Omega_2$ is nearly convex  and $\mbox{\ri}(\Omega_1\cap \Omega_2)=\ri\Omega_1\cap \ri\Omega_2$ (see~\cite[Theorem~2.1]{bkw2008} or~\cite[Corollary~4.8]{mmw2016}). Observe that
\begin{equation*}
F^{-1}(\Theta)=\mathcal{P}_1(\Omega_1\cap \Omega_2),
\end{equation*}  where $\mathcal{P}_1(x,y)=x$ for every $(x,y)\in\R^n\times\R^p$. Hence, by the just-cited results in~\cite{bkw2008} and~\cite{mmw2016}, we have
\begin{equation}\label{inequalities}
\mbox{\rm ri}\big( F^{-1}(\Theta)\big)=\mathcal{P}_1(\mbox{\ri}(\Omega_1\cap \Omega_2))=\mathcal{P}_1(\ri\Omega_1\cap \ri\Omega_2).
\end{equation}
Applying Theorem~\ref{Rthm}, we can assert that $(x, y)\in \ri\Omega_1\cap \ri\Omega_2$ if and only if $x\in \mbox{\rm ri}(\dom F)$  and $y\in \ri F(x)\cap \ri \Theta=F_0(x)\cap \ri\Theta$.  Since the inclusion $y\in F_0(x)\cap \ri\Theta$ holds for some $y\in\R^p$ if and only if $x\in F_0^{-1}(\ri\Theta)$,~\eqref{inequalities} implies~\eqref{R2}. $\h$

\medskip
Now, we use Theorem \ref{RINV} to derive  a significant result due to Bo{\c{t}}, Kassay  and Wanka~\cite{bkw2008}.

\begin{Corollary}\label{riinv}{\rm\bf (\cite[Theorem~2.2]{bkw2008})} Let $A\colon \R^n\to \R^p$ be a linear mapping and let $\Theta\subset\R^p$ be a nearly convex set.  If $A^{-1}(\ri \Theta)\neq\emptyset$, then
\begin{equation}\label{bkw_thm2.2}
\ri A^{-1}(\Theta)=A^{-1}(\ri \Theta).
\end{equation}
\end{Corollary}
{\bf Proof.} Let $F(x)=\{Ax\}$ for $x\in \R^n$. Then $\dom F=\R^n$ and $\rge F=F(\R^n)$ is an affine set, so $\mbox{\rm ri}(\rge F)=\rge F=F(\R^n)$. Choosing $\hat{x}\in A^{-1}(\ri \Theta)$, we see that  $\hat{y}=A(\hat{x})$ belongs to $\rge F\cap \ri\Theta$, so~\eqref{QC3} is satisfied in this setting. Observe that $F_0(x)=F(x)=\{Ax\}$ for all $x\in \R^n$. Thus,
\begin{equation*}
\ri A^{-1}(\Theta)=\ri F^{-1}(\Theta)=\mbox{\rm ri}(\dom F)\cap F_0^{-1}(\ri \Theta)=\R^n\cap A^{-1}(\ri\Theta)=A^{-1}(\ri\Theta),
\end{equation*} where the second equality is valid by~\eqref{R2}. This establishes~\eqref{bkw_thm2.2}. $\h$

\medskip
Given a set-valued mapping $F\colon \R^n\tto \R^p$ and a subset $\Omega\subset\R^n$, define the {\em restriction} of $F$ on $\Omega$ by
\begin{equation*}
	F_\Omega(x)=\begin{cases} F(x)&\mbox{\rm if }x\in \Omega,\\
		\emptyset&\mbox{\rm if }x\in \R^n\setminus \Omega.
	\end{cases}
\end{equation*}
The next theorem discusses the near convexity of $F_\Omega\colon \R^n\to \R^p$.

\begin{Theorem}\label{RNC} Let $F\colon \R^n\tto \R^p$ be a nearly convex set-valued mapping and let $\Omega$ be a nearly convex set in $\R^n$. Suppose that
	\begin{equation}\label{reguR}
		\mbox{\rm ri}(\dom F)\cap \ri \Omega\neq\emptyset.
	\end{equation}
	Then $F_\Omega$ is also nearly convex. In addition,
	\begin{equation}\label{restriction graph}
		\mbox{\rm ri}(\gph F_\Omega)=\big\{(x, y)\in \R^n\times \R^p\; \big |\; x\in \mbox{\rm ri}(\dom F)\cap \ri \Omega, \; y\in \ri F(x)\big\}.
	\end{equation}
\end{Theorem}
{\bf Proof.} Define the two nearly convex sets
\begin{equation*}
	\Omega_1=\gph F\; \mbox{\rm and }\; \Omega_2=\Omega\times \R^p.
\end{equation*}
Then we see that $\gph F_\Omega=\Omega_1\cap \Omega_2$. By \eqref{reguR} we can choose $\hat{x}\in \mbox{\rm ri}(\dom F)\cap \ri \Omega$. Since $F(\hat{x})$ is nonempty and nearly convex by Theorem \ref{Rthm}, the set $\ri F(\hat{x})$ is nonempty. Choosing $\hat{y}\in \ri F(\hat{x})$ and using Theorem \ref{Rthm} again show that $\ri\Omega_1\cap \ri\Omega_2\neq\emptyset$. Thus, $\gph F_\Omega=\Omega_1\cap \Omega_2$ is nearly convex. Representation \eqref{restriction graph} follows from Theorem \ref{Rthm}, Theorem \ref{T}(c), and the obvious fact that $\dom F_\Omega=(\dom F)\cap \Omega$. $\h$

Given a function $f\colon \R^n\to \oR$ and a subset $\Omega\subset\R^n$, define the {\em restriction} of $f$ on $\Omega$ by
\begin{equation*}
	f_\Omega(x)=\begin{cases} f(x)&\mbox{\rm if }x\in \Omega,\\
		\infty&\mbox{\rm if }x\in \R^n\setminus \Omega.
	\end{cases}
\end{equation*}
The following corollary is a direct consequence of Theorem~\ref{RNC}.

\begin{Corollary}\label{fr}  Let $f\colon \R^n\to \oR$ be a proper nearly convex function and let $\Omega$ be a nearly convex set in $\R^n$. Suppose that
	\begin{equation}\label{reguRf}
		\mbox{\rm ri}(\dom f)\cap \ri \Omega\neq\emptyset.
	\end{equation}
	Then $f_\Omega$ is also nearly convex. In addition,
	\begin{equation}\label{restriction epi}
		\mbox{\rm ri}(\epi f_\Omega)=\big\{(x, \lambda)\in \R^n\times \R\; \big |\; x\in \mbox{\rm ri}(\dom f)\cap \ri \Omega, \; f(x)<\lambda\big\}.
	\end{equation}
\end{Corollary}
{\bf Proof.} Consider the epigraphical mapping $F=E_f$ and see that $\gph F_\Omega=\epi f_\Omega$ and $\dom F_\Omega=\dom f_\Omega=(\dom f)\cap\Omega$. Theorem \ref{RNC} together with \eqref{reguRf} guarantees that $F_\Omega$ is nearly convex and thus the set $\gph F_\Omega=\epi f_\Omega$ is nearly convex. This implies the near convexity of $f_\Omega$. Representation \eqref{restriction epi} follows directly from \eqref{restriction graph}.  $\h$

\medskip
In the next theorem, we study the relative interior of the graph of the sum of two nearly convex set-valued mappings. As far as we know, there is no similar result in the literature.

\begin{Theorem} Let $F_1, F_2\colon \R^n\tto \R^p$ be nearly convex set-valued mappings. Suppose that
\begin{equation}\label{QC4}
\mbox{\rm ri}(\dom F_1)\cap \mbox{\rm ri}(\dom F_2)\neq\emptyset.
\end{equation}
Then we have
\begin{eqnarray}\label{risum}
\begin{array}{ll}
\mbox{\rm ri}(\mbox{\rm gph}(F_1+F_2))=\big\{(x, y)\in \R^n\times \R^p\; \big|\;& x\in \mbox{\rm ri}(\dom F_1)\cap \mbox{\rm ri}(\dom F_2),\\
& y\in \ri F_1(x)+\ri F_2(x)\big\}.
\end{array}
\end{eqnarray}
\end{Theorem}
{\bf Proof.} According to~\cite[Theorem~4.2]{NTY2023}, the qualification condition~\eqref{QC4} guarantees that $F_1+F_2$ is nearly convex. Since $\mbox{\rm dom}(F_1+F_2)=\dom F_1\cap \dom F_2$, under~\eqref{QC4} we have by~\cite[Theorem~2.1]{bkw2008} that
\begin{equation*}
\mbox{\rm ri}(\mbox{\rm dom}(F_1+F_2))=\mbox{\rm ri}(\dom F_1\cap \dom F_2)=\mbox{\rm ri}(\dom F_1)\cap \mbox{\rm ri}(\dom F_2).
\end{equation*}
By Theorem~\ref{Rthm}, for any $x\in \mbox{\rm ri}(\dom F_1)\cap \mbox{\rm ri}(\dom F_2)$ the sets $F_1(x)$ and  $F_2(x)$ are nearly convex. Thus, setting $T(y_1,y_2)=y_1+y_2$ for all $(y_1,y_2)\in \R^p\times\R^p$ and applying~\cite[Lemma~2.4]{bkw2008} yield
\begin{equation*}
\ri F(x)=\ri T\big[F_1(x)\times F_2(x)\big]=T\big[\big(\ri F_1(x)\big)\times \big(\ri F_2(x)\big)\big]=\ri F_1(x)+\ri F_2(x).
\end{equation*}
Therefore, representation~\eqref{risum} is a direct consequence of Theorem~\ref{Rthm}. $\h$

\medskip
Given a set-valued mapping $F\colon \R^n\tto \R^p$ and $G\colon \R^p\tto \R^q$, define
\begin{align*}
M_0(x, z)=\ri F(x) \cap \ri G^{-1}(z), \; (x, z)\in \R^n\times \R^p.
\end{align*}
Next, we study the relative interior of the graph of a composition of nearly convex set-valued mappings.
\begin{Theorem} Let $F\colon \R^n\tto \R^p$ and $G\colon \R^p\tto \R^q$ be nearly convex set-valued mappings. Suppose that
\begin{equation}\label{QC5}
\mbox{\rm ri}(\rge F)\cap \mbox{\rm ri}(\dom G)\neq\emptyset.
\end{equation}
Then we have
\begin{equation}\label{ricom}
\mbox{\rm ri}(\mbox{\rm gph}(G\circ F))=\big[\mbox{\rm ri}(\dom F)\times \mbox{\rm ri}(\rge G)\big]\cap \mbox{\rm dom}\,M_0.
\end{equation}
\end{Theorem}
{\bf Proof.} Let $\Omega_1=(\gph F)\times \R^q$ and let $\Omega_2=\R^n\times (\gph G)$. Define the projection mapping $\mathcal{P}_{1, 3}\colon \R^n\times \R^p\times \R^q\to \R^n\times \R^q$ by
\begin{equation*}
\mathcal{P}_{1, 3}(x, y, z)=(x, z), \ \; (x, y, z)\in \R^n\times \R^p\times \R^q.
\end{equation*}
Then we can show that
\begin{equation*}
\mbox{\rm gph}(G\circ F)=\mathcal{P}_{1,3}(\Omega_1\cap \Omega_2).
\end{equation*}
The  qualification condition \eqref{QC5} guarantees that
\begin{equation}\label{P13}
{\rm ri}(\mbox{\rm gph}(G\circ F))=\mathcal{P}_{1,3}(\ri\Omega_1\cap \ri\Omega_2).
\end{equation}
By Theorem \ref{Rthm} we have
\begin{align*}
&\mbox{\rm ri}(\Omega_1)=\{(x, y, z)\; |\; x\in \mbox{\rm ri}(\dom F), \; y\in \ri F(x)\},\\
&\mbox{\rm ri}(\Omega_2)=\{(x, y, z)\; |\; z\in \mbox{\rm ri}(\rge G), \; \; y\in \ri G^{-1}(z)\}.
\end{align*}
It follows that
\begin{equation*}
\ri\Omega_1\cap \ri\Omega_2=\big\{(x, y, z)\in \R^n\times \R^p\times \R^q\;\big |\; x\in \mbox{\rm ri}(\dom F),\; z\in \mbox{\rm ri}(\rge G), \;  \; y\in M_0(x, z)\big\}.
\end{equation*}
Representation~\eqref{ricom} now follows from \eqref{P13}. $\h$

\section{Near Convexity of Composite Set-Valued  Mappings}
\setcounter{equation}{0}

 In this section, we study the near convexity of a number of composite extended-real-valued  functions and set-valued mappings. Such functions and mappings are of independent interest, while being important for developing duality theory in the sequel.

 Given a set $\Theta\subset \R^n$ and set-valued mappings $F\colon \R^n\tto \R^p$, $G\colon \R^n\tto \R^q$, define
\begin{equation}\label{CS1}
\Phi(x,y)=\begin{cases} F(x)&\mbox{\rm if }x\in\Theta,\; y\in G(x),\\
\emptyset&\mbox{\rm otherwise}.
\end{cases}
\end{equation}
Then $\Phi\colon \R^n\times \R^q\tto\R^p$ is a composite set-valued mapping.

\begin{Theorem}\label{BOT11} Let $\Theta\subset\R^n$ be a nearly convex set, and let $F\colon \R^n\tto \R^p$, $G\colon \R^n\tto \R^q$ be  nearly convex set-valued mappings. Suppose that
\begin{equation}\label{RQC1}
\mbox{\rm ri}(\dom F)\cap\mbox{\rm ri}(\dom G) \cap \ri \Theta\neq\emptyset.
\end{equation}
Then the set-valued mapping $\Phi$ defined in \eqref{CS1} is nearly convex.
\end{Theorem}
{\bf Proof.} Define the set-valued mapping $H\colon \R^n\times \R^q\tto \R^p$ by
\begin{equation*}
 H(x, y)=F(x), \ \;  (x, y)\in \R^n\times \R^q.
\end{equation*}
 Then consider the set
\begin{equation*}
\Omega=\big\{(x, y)\in \R^n\times \R^q\; |\; x\in \Theta, \; y\in G(x) \big\},
\end{equation*}
and observe that $\Phi=H_\Omega$, i.e., $\Phi$ is the restriction of $H$ on $\Omega$.  We have \begin{equation}\label{ridomH}\mbox{\rm ri}(\dom H)=\mbox{\rm ri}(\dom F)\times \R^q.\end{equation} Using Theorem~\ref{T}(c) and Theorem~\ref{Rthm}, we can show that $\Omega$ is a nearly convex set with
\begin{equation}\label{riOmega4}
\ri\Omega=\big\{(x, y)\in \R^n\times \R^q\; \big |\; x\in \ri \Theta\cap \mbox{\rm ri}(\dom G), \; y\in \ri G(x)\big\}.
\end{equation}
Since $H$ is obviously nearly convex, it suffices to show that $\mbox{\rm ri}(\dom H)\cap \ri \Omega\neq \emptyset$ and then use Theorem~\ref{RNC}. By \eqref{RQC1}, there exists $\hat{x}\in \mbox{\rm ri}(\dom F)\cap\mbox{\rm ri}(\dom G)\cap  \ri \Theta$. Since $G$ is nearly convex, by Theorem~\ref{Rthm} the set $\ri G(\hat{x})$ is nonempty, so there exists $\hat{y}\in \ri G(\hat{x})$. Then, we see by \eqref{ridomH} and \eqref{riOmega4} that $(\hat{x}, \hat{y})\in \mbox{\rm ri}(\dom H)\cap \ri\Omega$. Therefore, the near convexity of $\Phi$ follows from Theorem~\ref{RNC}. $\h$

Theorem \ref{BOT11} can be applied to set-valued mappings defined by the intersection of a set constraint and an inequality constraint induced from a partially ordered relation
\begin{equation*}
y_1\leq_K y_2 \ \; \mbox{\rm if and only if }\; y_2-y_1\in K,
\end{equation*}
where $K$ is a nonempty convex cone. Namely, let $F\colon \R^n\tto \R^p$ be a set-valued mapping. For  $g\colon \R^n\to \R^q$ and $\Theta\subset\R^n$ define
\begin{equation}\label{CS2}
\Phi_1(x,y)=\begin{cases} F(x)&\mbox{\rm if }x\in \Theta,\; g(x)\leq_K y,\\
\emptyset&\mbox{\rm otherwise}.
\end{cases}
\end{equation}
The corollary below discusses the near convexity of $\Phi_1$.

\begin{Corollary}\label{BOT1} Let $\Theta\subset\R^n$ be a nearly convex set, $F\colon \R^n\tto \R^p$ a nearly convex set-valued mapping, and $g\colon \R^n\to \R^q$ a $K$-nearly convex function, where $K$ is a nonempty convex cone. Suppose that
\begin{equation}\label{RQC2}
\mbox{\rm ri}(\dom F)\cap \ri \Theta\neq\emptyset.
\end{equation}
Then the set-valued mapping $\Phi_1$ defined in \eqref{CS2} is nearly convex.
\end{Corollary}
{\bf Proof.} Define the set-valued mapping $G\colon \R^n\tto \R^q$ by
\begin{equation*}
G(x)=g(x)+K, \ \; x\in \R^n.
\end{equation*}
Then $G$ is nearly convex  because $g$ is $K$-nearly convex. Clearly, $\Phi_1$ has the form \eqref{CS1} with this particular set-valued mapping $G$. Since $\dom G=\R^n$, the qualification condition \eqref{RQC2} implies that
\begin{equation*}
\mbox{\rm ri}(\dom F)\cap\mbox{\rm ri}(\dom G) \cap \ri \Theta=\mbox{\rm ri}(\dom F)\cap \ri \Theta\neq\emptyset.
\end{equation*}
Therefore, the near convexity of $\Phi_1$ follows from Theorem~\ref{BOT11}. $\h$

Next, we discuss the near convexity for the restriction of a nearly convex function $f\colon \R^n\to \oR$ on a nearly convex constraint set. Given  a set-valued mapping $G\colon \R^n\tto \R^q$ and a set $\Theta\subset\R^n$,  define
\begin{equation}\label{P22}
\phi(x,y)=\begin{cases} f(x)&\mbox{\rm if }x\in \Theta,\; y\in G(x),\\
\infty&\mbox{\rm otherwise}.
\end{cases}
\end{equation}
The corollary below is another direct consequence of Theorem~\ref{BOT11}.

\begin{Corollary}\label{bot3}  Let $f\colon \R^n\to \oR$ be a proper nearly convex function, $G\colon \R^n\tto \R^q$ a nearly convex set-valued mapping, and $\Theta\subset\R^n$ a nearly convex set. Suppose that
\begin{equation}\label{RQCf1}
\mbox{\rm ri}(\dom f)\cap \mbox{\rm ri}(\dom G) \cap \ri \Theta\neq\emptyset.
\end{equation}
Then the function $\phi$ defined in \eqref{P22} is nearly convex.
\end{Corollary}
{\bf Proof.}  Consider the set-valued mapping $\Phi$ from \eqref{CS1}, where $F$ is the epigraphical mapping associated with $f$, i.e., $F=E_f$. Since $\dom F=\dom f$, we see that the qualification condition \eqref{RQC1} is satisfied under \eqref{RQCf1}. We deduce from Theorem~\ref{BOT11} that $\Phi$ is nearly convex. Thus, $\gph\Phi$ is a nearly convex set. It follows directly from the definition that $\gph \Phi=\epi \phi$. Therefore, $\phi$ is nearly convex since $\epi \phi$ is nearly convex.   $\h$

We will now discuss the restriction of a function $f\colon \mathbb{R}^n\rightarrow \overline{\mathbb{R}}$ on a constraint set given as the intersection of a geometric constraint and an inequality constraint from equation \eqref{CS2}. Given $g\colon \mathbb{R}^n\rightarrow\mathbb{R}^q$ and $\Theta\subset\mathbb{R}^n$, define
\begin{equation}\label{P2}
\phi_1(x,y)=\begin{cases} f(x)&\mbox{\rm if }x\in \Theta,\; g(x)\leq_K y,\\
\infty&\mbox{\rm otherwise}.
\end{cases}
\end{equation}

The following corollary presents another result from \cite{bkw2008} that directly follows from Corollary~\ref{bot3}.

\begin{Corollary}{\rm\bf (\cite[Theorem~4.1]{bkw2008})} Let $f\colon \R^n\to \oR$ be a proper nearly convex function, $g\colon \R^n\to \R^q$ a $K$-nearly convex function, and $\Theta\subset\R^n$ a nearly convex set. Suppose that
\begin{equation}\label{RQC3}
\mbox{\rm ri}(\dom f)\cap \ri \Theta\neq\emptyset.
\end{equation}
Then the function $\phi_1$ defined in \eqref{P2} is nearly convex.
\end{Corollary}
{\bf Proof.} It suffices to consider $G(x)=g(x)+K$, $x\in \R^n$, and then apply Corollary~\ref{bot3}.  $\h$

Next, we use the framework at the beginning of this section and define $\Psi\colon \R^n\times \R^n\times \R^q\tto \R^p$ by
\begin{equation}\label{P33}
\Psi(x,u,y)=\begin{cases} F(x+u)&\mbox{\rm if }x\in \Theta,\; y\in G(x),\\
\emptyset &\mbox{\rm otherwise}.
\end{cases}
\end{equation}

The following theorem demonstrates that the set-valued mapping $\Psi$ is nearly convex under the same qualification considered in Theorem~\ref{BOT11}.

\begin{Theorem}\label{BOT23} Let $F\colon \R^n\tto \R^p$, $G\colon \R^n\tto \R^q$ be nearly convex set-valued mappings, and let $\Theta\subset\R^n$ be a nearly convex set.  Suppose that the qualification condition~\eqref{RQC1} from Theorem~\ref{BOT11} is satisfied. Then the set-valued mapping $\Psi$ defined in \eqref{P33} is nearly convex.
\end{Theorem}
{\bf Proof.}  Let $H(x, u, y)=F(x+u)$ for $(x, u, y)\in \R^n\times \R^n\times \R^p$ and consider the set
\begin{equation*}
\Omega=\big\{(x, u, y)\in \R^n\times \R^n\times \R^p\; \big |\; x\in \Theta, \; y\in G(x)\big\}.
\end{equation*}
Since $\mbox{\rm ri}(\dom G) \cap \ri \Theta\neq\emptyset$ by~\eqref{RQC1}, we deduce from Theorem~\ref{T}(c) and Theorem~\ref{Rthm} that $\Omega$ is nearly convex with
\begin{equation*}
\ri\Omega=\big\{(x, u, y)\in \R^n\times \R^n\times \R^q\; \big |\; x\in\ri  \Theta, \; y\in \ri G(x)\big\}.
\end{equation*}
By the near convexity of the composition from \cite[Theorem~4.6]{NTY2023}, the set-valued mapping $H$ is nearly convex. Observe that $\Psi=H_\Omega$.  Let us now show that $\mbox{\rm ri}(\dom H)\cap \ri\Omega\neq\emptyset$. Indeed, we have
\begin{equation*}
\dom H=\big\{(x, u, y)\; \big |\; x+u\in \dom F, \; y\in \R^q\big\}.
\end{equation*}
Consider the linear mapping $A\colon \R^n\times \R^n\to \R^n$ given by $A(x, u)=x+u$ for $(x, u)\in \R^n\times \R^n$ and choose $\hat{x}\in \mbox{\rm ri}(\dom F)$. Then $(\hat{x}, 0)\in A^{-1}(\mbox{\rm ri}(\dom F){\color{blue})}$. Using Corollary~\ref{riinv}, we have
\begin{equation*}
\mbox{\rm ri}(\dom H)=\big\{(x, u, y)\; \big |\; x+u\in \mbox{\rm ri}(\dom F), \; y\in \R^q\big\}.
\end{equation*}
Now, choose $\hat{x}\in \mbox{\rm ri}(\dom F)\cap \mbox{\rm ri}(\dom G)\cap \ri \Theta$. Then choose $\hat{y}\in \ri G(\hat{x})$. By the representations of $\ri \Omega$ and $\mbox{\rm ri}(\dom H)$ we have $(\hat{x}, 0, \hat{y})\in \mbox{\rm ri}(\dom H)\cap \ri \Omega$. By Theorem~\ref{RNC} we see that $\Psi=H_\Omega$ is nearly convex. $\h$

In the same context as Corollary~\ref{BOT1}, we examine a specific instance of the set-valued mapping $\Psi$ from~\eqref{P33}. Define $\Psi_1\colon \mathbb{R}^n\times \mathbb{R}^n\times \mathbb{R}^q\rightrightarrows \mathbb{R}^p$ as follows:
\begin{equation}\label{P3}
\Psi_1(x,u,y)=\begin{cases} F(x+u)&\mbox{\rm if }x\in \Theta,\; g(x)\leq_K y,\\
\emptyset &\mbox{\rm otherwise}.
\end{cases}
\end{equation}

\begin{Corollary}\label{BOT2} Let $F\colon \R^n\tto \R^p$ be a nearly convex set-valued mapping, let $\Theta\subset\R^n$ be a nearly convex set, and let $g\colon \R^n\to \R^q$ be a $K$-nearly convex function, where $K$ is a nonempty convex cone. Suppose that the qualification condition \eqref{RQC2} from Corollary~\ref{BOT1} is satisfied. Then the set-valued mapping $\Psi_1$ defined in \eqref{P3} is nearly convex.
\end{Corollary}
{\bf Proof.}  Note that $\Psi_1$ is a particular case of the set-valued mapping $\Psi$ defined in~\eqref{P33}, where $G(x)=g(x)+K$ for $x\in \R^n$. Since $G$ is nearly convex with $\dom G=\R^n$ under the $K$-near convexity of $g$, the conclusion follows directly from Theorem~\ref{BOT23}. $\h$

We will now consider a new function that will be useful for our subsequent discussion on {\em  Fenchel-Lagrange duality}. This function is similar to the function $\phi_1$ defined in~\eqref{P22}, but includes an additional auxiliary variable $u$:

\begin{equation}\label{P44}
\psi(x,u,y)=\begin{cases} f(x+u)&\mbox{\rm if }x\in \Theta, \; y\in G(x),\\
\infty&\mbox{\rm otherwise}.
\end{cases}
\end{equation}

\begin{Corollary}\label{bot33}  Let $f\colon \R^n\to \oR$ be a nearly convex function, let $G\colon \R^n\tto \R^q$ be a nearly convex set-valued mapping, and let  $\Theta\subset\R^n$ be a nearly convex set.  Suppose that
\begin{equation*}
\mbox{\rm ri}(\dom f)\cap \mbox{\rm ri}(\dom G)\cap \ri \Theta\neq\emptyset.
\end{equation*}
Then the function $\psi$ defined in \eqref{P44} is nearly convex.
\end{Corollary}
{\bf Proof.} We use the same arguments as in the proof of Corollary~\ref{bot3}. It suffices to consider $F=E_f$ and then apply Theorem~\ref{BOT23}.  $\h$

Next, we consider a particular case of the function $\phi$ from \eqref{P44}, where $G$ is a {\em $K$-epigraphical mapping}.
\begin{equation}\label{P4}
\psi_1(x,u,y)=\begin{cases} f(x+u)&\mbox{\rm if }x\in \Theta,\; g(x)\leq_K y,\\
\infty&\mbox{\rm otherwise}.
\end{cases}
\end{equation}
The proof of the next corollary is a straightforward application of Corollary~\ref{bot33}.
\begin{Corollary}{\rm\bf (\cite[Theorem~4.3]{bkw2008})} Let $f\colon \R^n\to \oR$ be a nearly convex function, $\Theta\subset\R^n$ a nearly convex set, and $g\colon \R^n\to \R^p$ a $K$-nearly convex function, where $K$ is a nonempty convex cone.  Suppose that \eqref{RQC3} is satisfied. Then the function $\psi_1$ defined in \eqref{P4} is nearly convex.
\end{Corollary}
{\bf Proof.} It suffices to consider $G(x)=g(x)+K$ for $x\in \R^n$,  and then use Corollary~\ref{bot33}.  $\h$

In the following theorem, we examine the near convexity of a composite set-valued mapping with two variables. This result differs from the near convexity of sums and compositions of set-valued mappings, as presented in \cite[Theorem~4.2]{NTY2023} and \cite[Theorem~4.6]{NTY2023} because it does not require any qualification conditions, although it is a consequence of that theorems.

\begin{Theorem}\label{thm_sum1} Let $F\colon \R^n\tto \R^q$ and $G\colon \R^p\tto \R^q$  be proper nearly convex set-valued mappings. Define
	\begin{equation*}
		\Phi(x, y)=F(x)+G(Ax+y), \ \; (x, y)\in \R^n\times \R^p,
	\end{equation*}
	where $A\in \R^{p\times n}$. Then $\Phi\colon \R^n\times \R^p\tto \R^q$ is a nearly convex set-valued mapping.
\end{Theorem}
{\bf Proof.} Let $\Phi_1(x, y)=F(x)$ and let $\Phi_2(x, y)=G(Ax+y)$ for $(x, y)\in \R^n\times \R^p$. It follows directly from the definition that $\Phi_1$ is nearly convex with  $\dom \Phi_1=(\dom F)\times \R^p$. Define the linear mapping
\begin{equation*}
	T(x, y)=Ax+y, \ \; (x, y)\in \R^n\times \R^p.
\end{equation*}
Observe that  $\Phi_2(x, y)=G(T(x, y))$ for $x, y\in \R^n\times \R^p$ and that $T^{-1}(\mbox{\rm ri}(\dom G))\neq\emptyset$. Indeed, let $\hat{u}\in \mbox{\rm ri}(\dom G)$, which is a nonempty set due to the near convexity of $\dom G$. Then $(0, \hat{u})\in T^{-1}(\mbox{\rm ri}(\dom G))$. By \cite[Theorem~4.6]{NTY2023}, the set-valued mapping $\Phi_2$ is nearly convex with $\dom \Phi_2=T^{-1}(\dom G)$. To show that $\Phi$ is nearly convex, it suffices to show that $\mbox{\rm ri}(\dom \Phi_1)\cap \mbox{\rm ri}(\dom\Phi_2)\neq\emptyset$ and then use \cite[Theorem~4.2]{NTY2023}. Obviously, $$\mbox{\rm ri}(\dom \Phi_1)=\mbox{\rm ri}(\dom F)\times \R^p.$$ To compute $\mbox{\rm ri}(\dom \Phi_2)$, we use  Corollary~\ref{riinv} and get
\begin{equation*}
	\mbox{\rm ri}(\dom \Phi_2)=\mbox{\rm ri}(T^{-1}(\dom G))=T^{-1}(\mbox{\rm ri}(\dom G))=\{(x, y)\in \R^n\times \R^p \; |\; Ax+y\in \mbox{\rm ri}(\dom G)\}.
\end{equation*}
Now choose $\hat{x}\in \mbox{\rm ri}(\dom F)$ and let $\hat{y}=\hat{u}-A\hat x$. It is obvious that $$(\hat{x}, \hat{y})\in \mbox{\rm ri}(\dom \Phi_1)\cap \mbox{\rm ri}(\dom\Phi_2),$$ which justifies the near convexity of $\Phi$ by~\cite[Theorem~4.2]{NTY2023}. $\h$

As a corollary, we obtain another result by Bo{\c{t}}, Kassay and Wanka.

\begin{Corollary}{\rm\bf (\cite[Theorem~4.5]{bkw2008})} Let $f\colon \R^n\to \oR$ and  $g\colon \R^p\to \oR$ be proper nearly convex functions. Define
	\begin{equation*}
		\phi(x, y)=f(x)+g(Ax+y), \ \; (x, y)\in \R^n\times \R^p,
	\end{equation*}
	where $A\in \R^{p\times n}$. Then $\phi$ is nearly convex.
\end{Corollary}
{\bf Proof.} Consider the epigraphical mappings $F=E_f\colon \R^n\tto\R$ and $G=E_g\colon \R^p\tto \R$. Then, by Theorem~\ref{thm_sum1}, the set-valued mapping defined by
\begin{equation*}
	\Phi(x, y)=E_f(x)+E_g(Ax+y)=[f(x)+g(Ax+y), \infty), \ \; (x, y)\in \R^n\times \R^p,
\end{equation*}
is nearly convex. Thus, $\gph \Phi=\epi \phi$ is nearly convex, so $\phi$ is nearly convex. $\h$

\section{Subdifferentiation of Optimal Value Functions under Near Convexity}
\setcounter{equation}{0}

This section focuses on studying a class of extended-real-valued functions known as the {\em  optimal value function} under near convexity. We demonstrate that when an objective function is nearly convex and a set-valued mapping is nearly convex, the resulting optimal value function is also nearly convex, provided that a relative interior qualification condition is satisfied.

Given a function $f\colon \R^n\times \R^p\to \oR$ and a set-valued mapping $F\colon \R^n\tto \R^p$, we define the  optimal value function
\begin{equation}\label{mf}
\mu(x)=\inf\{f(x, y)\; |\; y \in F(x)\}, \ \; x\in \R^n.
\end{equation}
This class of functions plays a crucial role in many aspects of variational and convex analysis. The {\em strict epigraph} of an extended-real-valued function $\ph\colon \R^n\to \oR$ {\color{blue} is} defined by
\begin{equation*}
\epi_s\, \ph=\big\{(x, \gamma)\in \R^n\times \R\; \big|\; \ph(x)<\gamma\big\}.
\end{equation*}
The lemma below provides a useful relationship between the strict epigraph and the graph of a function.
\begin{Lemma}\label{sepi} Let $\ph\colon \R^n\to \oR$ be a function. Then we have the equality
\begin{equation}\label{EIQ}
\Bar{\epi_s\, \ph}=\Bar{\epi \ph}.
\end{equation}
\end{Lemma}
{\bf Proof.} Since $\epi_s\, \ph\subset \epi \ph$, we see that $\Bar{\epi_s\, \ph}\subset\Bar{\epi \ph}$. To justify the reverse inclusion, it suffices to show that
\begin{equation}\label{RI}
\epi \ph\subset \Bar{\epi_s\, \ph}.
\end{equation}
For any $(x, \lambda)\in \epi \ph$ we have $(x, \lambda)\in \R^n\times \R$ and $\ph(x)\leq \lambda$. First, consider the case where $-\infty <\ph(x)$, which implies that $\ph(x)\in \R$. Letting $\lambda_k=\ph(x)+1/k$ for $k\in \N$, we have $(x, \lambda_k)\in \epi_s \, \ph$ and $(x, \lambda_k)\to (x, \lambda)$ as $k\to \infty$. Thus, $(x, \lambda)\in \Bar{\epi_s\, \ph}$. Now, consider the case where $\ph(x)=-\infty$. In this case it is obvious that $-\infty=\ph(x)<\lambda$. Then $(x, \lambda)\in \epi_s\, \ph\subset \Bar{\epi_s\, \ph}$. Therefore, \eqref{RI} is satisfied, which justifies equality \eqref{EIQ} and completes the proof. $\h$

 Our result on the near convexity of the optimal value function is formulated as follows.

\begin{Theorem}\label{mfnc} Let $f\colon \R^n\times \R^p\to \oR$ is a proper nearly convex function, and let $F\colon \R^n\tto \R^p$ be a nearly convex set-valued mapping. Suppose that
\begin{equation}\label{QC1}
\mbox{\rm ri}(\dom f)\cap \mbox{\rm ri}(\gph F)\neq\emptyset.
\end{equation}
Then the optimal value function $\mu$ given in \eqref{mf} is nearly convex.
\end{Theorem}
{\bf Proof.} Define the nearly convex sets
\begin{equation}\label{om12}
\Omega_1=\epi f,\ \; \Omega_2=(\gph F)\times \R.
\end{equation}
We first claim that
\begin{equation}\label{I1}
\epi_s\, \mu\subset\mathcal{P}_{1,3}(\Omega_1\cap \Omega_2)\subset \epi\mu,
\end{equation}
where $\mathcal{P}_{1,3}\colon \R^n\times \R^p\times \R\to \R^n\times \R$ is the projection mapping defined by
\begin{equation*}
\mathcal{P}_{1, 3}(x, y, \lambda)=(x, \lambda), \ \; (x, y, \lambda)\in \R^n\times \R^p\times \R.
\end{equation*}
Indeed, fix any $(x_0, \gamma_0)\in \epi_s\, \mu$ and get $\mu(x_0)=\inf_{y\in F(x_0)}f(x_0, y)<\gamma_0$. Thus, there exists $y_0\in F(x_0)$ such that
\begin{equation*}
y_0\in F(x_0)\ \; \mbox{\rm and }\; f(x_0, y_0)<\gamma_0.
\end{equation*}
Then $(x_0, y_0, \gamma_0)\in \Omega_1\cap \Omega_2$. This justifies the first inclusion in \eqref{I1}. To prove the second inclusion, fix any $(x_0, \gamma_0)\in \mathcal{P}_{1,3}(\Omega_1\cap \Omega_2)$. Then there exists $y_0\in \R^p$ such that
\begin{equation*}
(x_0, y_0, \gamma_0)\in \Omega_1\cap \Omega_2.
\end{equation*}
By the definitions of the sets $\Omega_1$ and $\Omega_2$ we have $f(x_0, y_0)\leq \gamma_0$ and also $y_0\in F(x_0)$. Thus
\begin{equation*}
\mu(x_0)=\inf_{y\in F(x_0)}f(x_0, y)\leq f(x_0, y_0)\leq \gamma_0,
\end{equation*}
which shows that $(x_0, \gamma_0)\in \epi \mu$ and completes the proof of \eqref{I1}. By \eqref{QC1} we can choose
\begin{equation*}
(\hat{x}, \hat{y})\in \mbox{\rm ri}(\dom f)\cap \mbox{\rm ri}(\gph F).
\end{equation*}
Then choosing $\hat{\gamma}>f(\hat{x}, \hat{y})$, we see that $(\hat{x}, \hat{y}, \hat{\gamma})\in \ri \Omega_1\cap \ri \Omega_2$, so $\ri \Omega_1\cap \ri\Omega_2\neq\emptyset$. Since both $\Omega_1$ and $\Omega_2$ are nearly convex, we can deduce from Theorem~\ref{T}(c) that $\Omega_1\cap \Omega$ is nearly convex. Thus,  $\mathcal{P}_{1,3}(\Omega_1\cap \Omega_2)$ is nearly convex by Theorem~\ref{T}(d).

Letting $C=\Bar{\mathcal{P}_{1,3}(\Omega_1\cap \Omega_2)}$, it follows from \eqref{I1} and Lemma~\ref{sepi} that
\begin{equation*}
\Bar{\epi_s\, \mu}=\Bar{\epi \mu}\subset C\subset \Bar{\epi \mu}.
\end{equation*}
Thus, $\Bar{\epi \mu}=C$, which is a closed convex set by the near convexity of $\mathcal{P}_{1,3}(\Omega_1\cap \Omega_2)$. Then
\begin{equation}\label{C1}
\epi \mu\subset\Bar{\epi \mu}=\Bar{\epi_s}\,\mu= C=\Bar{C}.
\end{equation}
Since $\mathcal{P}_{1,3}(\Omega_1\cap \Omega_2)$ is nearly convex, it follows from Theorem~\ref{T} that
\begin{align*}
\ri C&=\mbox{\rm ri}\big(\Bar{\mathcal{P}_{1,3}(\Omega_1\cap \Omega_2)}\big)=\mbox{\rm ri}\big(\mathcal{P}_{1,3}(\Omega_1\cap \Omega_2)\big)\\
&=\mathcal{P}_{1,3}\big(\mbox{\rm ri}(\Omega_1\cap \Omega_2)\big)=\mathcal{P}_{1, 3}\big(\ri \Omega_1\cap \ri\Omega_2\big).
\end{align*}
If $(x, y, \lambda)\in \ri \Omega_1\cap \ri\Omega_2$, then by \cite[Proposition~3.7]{NTY2023} we have $f(x, y)<\lambda$ with $y\in F(x)$, so $\mu(x)<\lambda$. Thus,
\begin{equation}\label{C2}
\ri C= \mathcal{P}_{1, 3}(\ri \Omega_1\cap \ri\Omega_2)\subset \epi_s\, \mu\subset \epi \mu.
\end{equation}
Combining \eqref{C1} and \eqref{C2} yields the near convexity of $\epi \mu$, which gives us the near convexity of $\mu$. $\h$

Next, we will obtain a formula for the subdifferential of the optimal value function \eqref{mf} defined by a nearly convex objective function and a nearly convex constraint mapping. This formula extends the known result for the convex case with a significantly different proof based on a geometric approach. We need to review some concepts of generalized differentiation.

Given a nearly convex set $\Omega$ in $\R^n$ with $\ox\in \Omega$, define the {\em normal cone} to $\Omega$ at $\ox$ by
\begin{equation*}
N(\ox; \Omega)=\big\{v\in \R^n\; |\; \la v, x-\ox\ra\leq 0\ \, \mbox{\rm for all }\, x\in \Omega\big\}.
\end{equation*}
 Let $F\colon \R^n\tto \R^p$ be a nearly convex set-valued mapping and let $(\ox,\oy)\in \gph F$. The {\em coderivative} of $F$ at $(\ox,\oy)$ is the set-valued mapping $D^*F(\ox,\oy)\colon \R^p\tto \R^n$ with the values
	\begin{equation*}\label{cod}
		D^*F(\ox,\oy)(v)=\big\{u\in \R^n\;\big|\;(u,-v)\in N\big((\ox,\oy);\gph F\big)\big\},\ \, v\in \R^p.
	\end{equation*}
Let $\ph\colon \R^n\to \oR$ be a nearly convex function. We define the {\em subdifferential} of $\ph$ at $\ox\in \R^n$ with $\ph(\ox)\in \R$ by
\begin{equation*}
\partial \ph(\ox)=\big \{v\in \R^n\; \big|\; \la v, x-\ox\ra\leq \ph(x)-\ph(\ox)\ \, \mbox{\rm for all }\, x\in \R^n \big\}.
\end{equation*}
We can show that if $\ph\colon \R^n \to \oR$ is a  nearly convex function with $\ph(\ox)\in \R$, then
\begin{equation*}
D^*E_\ph(\ox, \ph(\ox))(1)=\partial \ph(\ox),
\end{equation*}
where $E_\ph$ is defined in \eqref{Emapping}.

In connection with the parametric minimization problem in~\eqref{mf}, we consider the solution map $S\colon \R^n\tto \R^p$ with
\begin{equation}\label{sol}
S(x):=\big\{y\in F(x)\; \big |\; f(x, y)=\mu(x)\big\}.
\end{equation}
We are now ready to establish the second main result of this section.
\begin{Theorem} Consider the optimal value function~\eqref{mf} in which $f$ is a proper nearly convex function and $F$ is a nearly convex set-valued mapping. Suppose that the qualification condition \eqref{QC1} is satisfied.  Take any $\ox\in \R^n$ with $\mu(\ox)\in \R$ and assume that $S(\ox)\neq\emptyset$.  Then for any $\oy\in S(\ox)$ we have
\begin{equation}\label{mfr}
\partial \mu(\ox)=\bigcup_{(u, v)\in \partial f(\ox, \oy)} \big[u+D^*F(\ox, \oy)(v)\big].
\end{equation} 
\end{Theorem}
{\bf Proof.} Fix any $\oy\in S(\ox)$ and take any $w\in \partial \mu(\ox)$. Then $$(w, 0, -1)\in N((\ox, \oy, \mu(\ox)); \Omega_1\cap \Omega_2),$$ where the sets $\Omega_1$ and $\Omega_2$ are given in \eqref{om12}. Indeed, for any $(x, y, \lambda)\in \Omega_1\cap \Omega_2$  we have
\begin{equation*}
f(x, y)\leq \lambda\ \; \mbox{\rm and }\; y\in F(x),
\end{equation*}
which implies that $\mu(x)\leq f(x, y)\leq \lambda$. Thus,
\begin{equation*}
\la w, x-\ox\ra+\la 0, y-\oy\ra+(-1)(\lambda-\mu(\ox))\leq \la w,  x-\ox\ra-(\mu(x)-\mu(\ox))\leq 0.
\end{equation*}
Since $\ri\Omega_1\cap \ri\Omega_2\neq\emptyset$ as in the proof of Theorem~\ref{mfnc}, we can use \cite[Theorem~5.1]{NTY2023} and get
\begin{align*}
(w, 0, -1)&\in N((\ox, \oy, \mu(\ox)); \Omega_1)+N((\ox, \oy, \mu(\ox)); \Omega_2)\\
&=N((\ox, \oy, f(\ox, \oy)); \epi f)+N((\ox, \oy, \mu(\ox)); \Omega_2).
\end{align*}
Thus, we have the representation
\begin{equation*}
(w, 0, -1)=(u, v, -1)+(u_1, -v, 0),
\end{equation*}
where $(u, v, -1)\in N((\ox, \oy, f(\ox, \oy)); \epi f)$ and $(u_1, -v)\in N((\ox, \oy); \gph F)$. It follows that $(u, v)\in \partial f(\ox, \oy)$ and $u_1\in D^*F(\ox, \oy)(v)$. Then
\begin{equation*}
w=u+u_1\in u+D^*F(\ox, \oy)(v)\subset \bigcup_{(u, v)\in \partial f(\ox, \oy)} \big[u+D^*F(\ox, \oy)(v)\big].
\end{equation*}
This justifies the inclusion $\subset$ in \eqref{mfr}. The proof of the reverse inclusion in~\eqref{mfr} is straightforward. $\h$

\section{Fenchel Conjugate Calculus under Near Convexity}
\setcounter{equation}{0}

This section aims to establish a formula for computing the Fenchel conjugate of the optimal value function \eqref{mf} when dealing with nearly convex data. This formula will then be utilized to derive several calculus rules for Fenchel conjugate of nearly convex functions.

Given two proper functions $g, h\colon \R^n\to \oR$, define the {\em infimal convolution} of these functions by
\begin{equation*}
(g\s h)(x)=\inf\{g(x_1)+h(x_2)\; |\; x_1+x_2=x\}, \ \; x\in \R^n.
\end{equation*}
For a nonempty subset $\Omega$ of $\R^n$, recall that the {\em support function} associated with $\Omega$ is defined by
\begin{equation*}
\sigma_{\small \Omega}(v)=\sup\{\la v, x\ra\; |\; x\in \Omega\}, \ \; v\in \R^n.
\end{equation*}
It follows from the definition that $\sigma_{\small\Omega}\colon \R^n\to \oR$ is a proper convex function regardless the nonconvexity of $\Omega$. In addition,
\begin{equation*}
\sigma_{\small \Omega}=\sigma_{\small \co\Omega}=\sigma_{\small \Bar\Omega}=\sigma_{\small \Bar{\co}\Omega}.
\end{equation*}
The theorem below allows us to represent the support function associated with the intersection of two nearly convex sets in terms of the infimal convolution of the support functions associated with those sets.

\begin{Theorem}\label{Sintersection} Let $\Omega_1$ and $\Omega_2$ be two nearly convex sets in $\R^n$. Assume that
\begin{equation*}
\ri \Omega_1\cap \ri\Omega_2\neq\emptyset.
\end{equation*}
Then we have the equality
\begin{equation*}
\sigma_{\small{\Omega_1\cap \Omega_2}}=\sigma_{\small{\Omega_1}}\s \sigma_{\small{\Omega_2}}.
\end{equation*}
In addition, for any $v\in \R^n$ such that $\sigma_{\small{\Omega_1\cap \Omega_2}}(v)\in \R$, there exist $v_1, v_2\in \R^n$ such that $v=v_1+v_2$ and
\begin{equation*}
\sigma_{\small{\Omega_1\cap \Omega_2}}(v)=\sigma_{\small{\Omega_1}}(v_1)+\sigma_{\small{\Omega_2}}(v_2).
\end{equation*}
\end{Theorem}
{\bf Proof.} Choose two convex sets $C_1$ and $C_2$ in $\R^n$ such that
\begin{equation*}
C_1\subset \Omega_1\subset\Bar{C_1}\; \mbox{\rm and }C_2\subset \Omega_2\subset\Bar{C_2}.
\end{equation*}
Then we have the equalities
\begin{equation*}
\ri C_i=\ri \Omega_i\; \mbox{\rm and }\Bar{C_i}=\Bar{\Omega_i}, \ \; i=1, 2.
\end{equation*}
It follows that
\begin{equation*}
\ri \Bar{C_1}\cap \ri \Bar{C_2}=\ri C_1\cap \ri C_2=\ri\Omega_1\cap \ri \Omega_2\neq\emptyset.
\end{equation*}
Then we can show that
\begin{equation}\label{clin}
\Bar{\Omega_1\cap \Omega_2}=\Bar{\Omega_1}\cap \Bar{\Omega_2}.
\end{equation}
Indeed, the inclusion $\subset$ in \eqref{clin} is obvious. By \cite[Theorem~2.26{\bf (i)}]{bmn} we have
\begin{equation*}
\Bar{\Omega_1}\cap \Bar{\Omega_2}=\Bar{C_1}\cap \Bar{C_2}=\Bar{C_1\cap C_2}\subset\Bar{\Omega_1\cap \Omega_2},
\end{equation*}
which justifies the reverse inclusion in \eqref{clin}. Using \eqref{clin} and \cite[Theorem~4.16]{bmn}, we have
\begin{align*}
\sigma_{\small{\Omega_1\cap \Omega_2}}&=\sigma_{\small{\Bar{\Omega_1\cap \Omega_2}}}=\sigma_{\small{\Bar{\Omega_1}\cap\Bar{\Omega_2}}}\\
&=\sigma_{\small{\Bar{C_1}\cap\Bar{C_2}}}=\sigma_{\small{\Bar{C_1}}}\s \sigma_{\small{\Bar{C_2}}}=\sigma_{\small{\Omega_1}}\s \sigma_{\small{\Omega_2}}.
\end{align*}
The last statement of the theorem can be justified similarly by using a related result for the convex case from \cite[Theorem~4.16]{bmn}. $\h$

Given a function $f\colon \R^n\to \oR$,  recall that  the {\em Fenchel conjugate} of $f$ is  given by
\begin{equation*}
f^*(w)=\sup\big\{\la w, x\ra-f(x)\; |\; x\in \R^n\big\}, \ \; w\in \R^n.
\end{equation*}

To continue, we introduce a new notion called the {\em Fenchel conjugate} of a set-valued mapping $F\colon \R^n\tto \R^p$. We define $F^*\colon \R^n\times \R^p\tto \oR$ by
\begin{equation*}
F^*(u, v)= \sigma_{\small\gph F}(u, v), \ \; (u, v)\in \R^n\times \R^p.
\end{equation*}
The proposition below allows us to represent the Fenchel conjugate of a function $f\colon \R^n\to\oR$ by the Fenchel conjugate of the epigraphical mapping.

\begin{Proposition}\label{fF} Let $f\colon \R^n\to \oR$ be  function. Then
\begin{equation*}
f^*(w)=E^*_f(w, -1),\ \; v\in \R^n,
\end{equation*}
where $E_f$ is the epigraphical mapping associated with $f$ defined in \eqref{Emapping}.
\end{Proposition}
{\bf Proof.} It follows from the definition that for any $w\in \R^n$ we have
\begin{equation*}
f^*(w)=\sigma_{\small\epi f}(w, -1)=\sigma_{\small\gph E_f}(w, -1)=(E_f)^*(w, -1),
\end{equation*}
which completes the proof. $\h$

The next theorem discusses the Fenchel conjugate of the optimal value function.
\begin{Theorem}\label{FC1} Consider the optimal value function \eqref{mf} in which $f$ is a proper nearly convex function and $F$ is a nearly convex set-valued mapping. Suppose that the qualification condition \eqref{QC1}. Then we have the representation
\begin{equation}\label{opF}
\mu^*(w)=(f^*\s F^*)(w, 0)\ \; \mbox{\rm for all }\; w\in \R^n.
\end{equation}
In addition, if $\mu^*(w)\in \R$, then there exist $w_1, w_2\in \R^n$ and $v\in \R^p$ such that
\begin{equation*}
\mu^*(w)=f^*(w_1, v)+F^*(w_2, -v).
\end{equation*}
\end{Theorem}
{\bf Proof.} Take any $w_1, w_2\in \R^n$ such that $w_1+w_2=w$ and take any $v\in \R^p$. For  $x, y\in \R^n$, we have the relationships
\begin{align*}
f^*(w_1, v)+F^*(w_2, -v)&=\sup\{\la w_1, x\ra +\la v, y\ra-f(x, y)\; |\; (x, y)\in \R^n\times \R^p\}\\
& \quad +\sup\{\la w_2, x\ra-\la v, y\ra\; |\; (x, y)\in \gph F\}\\\
&\geq w_1, x\ra +\la v, y\ra-f(x, y)+\la w_2, x\ra-\la v, y\ra\\
&=\la w, x\ra-f(x, y)\ \; \mbox{\rm whenever }y\in F(x).
\end{align*}
Taking the infimum with respect to all $y\in F(x)$ gives us
\begin{equation*}
f^*(w_1, v)+F^*(w_2, -v)\geq \la w, x\ra-\mu(x).
\end{equation*}
It follows by taking the supremum with respect to all $x\in \R^n$ that
\begin{equation*}
f^*(w_1, v)+F^*(w_2, -v)\geq \sup\big\{\la w, x\ra-\mu(x)\; \big |\; x\in \R^n\big\}=\mu^*(w).
\end{equation*}
Next, by the arbitrary choices of $w_1, w_2\in \R^n$ with $w_1+w_2=w$ and $v\in \R^p$ we have
\begin{equation}\label{opF1}
(f^*\s F^*)(w, 0)\geq \mu^*(w).
\end{equation}
Now, we will show that
\begin{equation}\label{opF11}
(f^*\s F^*)(w, 0)\leq \mu^*(w).
\end{equation}
If $\mu^*(w)=\infty$, then inequality \eqref{opF11} holds as an equality because $(f^*\s F^*)(w, 0)=\infty$ by \eqref{opF1}. The qualification condition \eqref{QC1} guarantees that $\mu^*(w)>-\infty$. Thus, we now only need to consider the case where $\mu^*(w)\in \R$. In this case, by Proposition~\ref{fF}  we see that
\begin{equation*}
\mu^*(w)=(E_\mu)^*(w, -1)=\sigma_{\Omega_1\cap\Omega_2}(w, 0, -1),
\end{equation*}
where $\Omega_1$ and $\Omega_2$ are nearly convex sets given in \eqref{om12}. Since $\ri\Omega_1\cap \ri \Omega_2\neq\emptyset$, we can employ Theorem~\ref{Sintersection} and find $w_1, w_2\in \R^n$ and $v\in \R^p$ such that
\begin{equation*}
\mu^*(w)=\sigma_{\Omega_1\cap\Omega_2}(w, 0, -1)=\sigma_{\Omega_1}(w_1, -v, \gamma_1)+\sigma_{\Omega_2}(w_2, -v, \gamma_2),
\end{equation*}
where $w_1+w_2=w$ and $\gamma_1+\gamma_2=-1$. Then $\gamma_2=0$ since $\sigma_{\Omega_2}(w_2, -v, \gamma_2)=\infty$ otherwise by the structure of $\Omega_2$. Thus, $\gamma_1=-1$, and so
\begin{align*}
\mu^*(w)&=\sigma_{\Omega_1\cap\Omega_2}(w, 0, -1)=\sigma_{\Omega_1}(w_1, v, -1)+\sigma_{\gph F}(w_2, -v)\\
&=f^*(w_1, v)+F^*(w_2, -v)\geq (f^*\s F^*)(w, 0).
\end{align*}
Combining this with \eqref{opF1} completes the proof.  Note that the proof of the last statement in the proof is included in our argument above.  $\h$

The next corollary considers a specific case of the optimal value function \eqref{mf} in which the set-valued mapping $F$ does not involve, i.e., the values of $F$ are the entire space $\R^p$.

\begin{Corollary}\label{FC1.1} Let $f\colon \R^n\times \R^p\to \oR$ be a proper nearly convex function. Define
\begin{equation*}\label{mf1.1}
\mu(x)=\inf\{f(x, y)\; |\; y\in \R^p\}, \ \; x\in \R^n.
\end{equation*}
Then we have the equality
\begin{equation}\label{mf1.1}
\mu^*(v)=f^*(v, 0)\; \mbox{\rm for all }v\in \R^n.
\end{equation}
\end{Corollary}
{\bf Proof.} Consider the constant set-valued mapping $F\colon \R^n\tto \R^p$ given by $F(x)=\R^p$ for all $x\in \R^n$ in the setting of Theorem~\ref{FC1}. Then $\gph F=\R^n\times \R^p$ and thus
\begin{equation*}
F^*(v, w)=\sup\{\la v, x\ra+\la w, y\ra\; |\; (x, y)\in \R^n\times \R^p\}=\begin{cases} 0\; &\mbox{\rm if }(v, w)=(0, 0),\\
\infty \; &\mbox{\rm otherwise}.
\end{cases}
\end{equation*}
Using Theorem~\ref{FC1}, we see that
\begin{equation*}
\mu^*(v)=(f^*\s F^*)(v, 0)=\inf\{f^*(v_1, -w)+F^*(v_2, w)\; |\; v_1+v_2=v, \; w\in \R^p\}=f^*(v, 0).
\end{equation*}
Indeed, since $F^*(v_2, w)=\infty$ if $v_2\neq 0$ or $w\neq 0$. Thus, we are only concerned with the case where $v_2=0$ and $w=0$. Then $v_1=v$ and $-w=0$. The qualification condition required in Theorem~\ref{FC1} is also satisfied in this case. $\h$

\begin{Remark}\label{FC1.12}{\rm It can be shown by a direct calculation that \eqref{mf1.1} holds for an arbitrary function $f\colon \R^n\to \oR$.}
\end{Remark}

The formula for the Fenchel conjugate of the optimal value function, as presented in Theorem~\ref{FC1}, has a broad applicability that encompasses a wide range of operations on functions and set-valued mappings. Using this formula, we now derive a Fenchel conjugate sum rule for nearly convex functions.
\begin{Theorem}\label{FC2}  Let $f_1, f_2\colon \R^n\to \oR$ be proper nearly convex functions. Suppose that
\begin{equation*}
\mbox{\rm ri}(\dom f_1)\cap \mbox{\rm ri}(\dom f_2)\neq\emptyset.
\end{equation*}
Then we have the equality
\begin{equation}\label{sumf}
(f_1+f_2)^*(w)=(f^*_1\s f^*_2)(w), \ \; w\in \R^n.
\end{equation}
In addition, if $(f_1+f_2)^*(w)\in \R$, then there exist $w_1, w_2\in \R^n$ such that $w=w_1+w_2$ and $(f_1+f_2)^*(w)=f^*_1(w_1)+f^*_2(w_2)$.
\end{Theorem}
{\bf Proof.} Consider the function $f\colon \R^n\times \R\to \oR$ given by $f(x, y)=f_1(x)+y$  and let $F(x)=E_{f_2}(x)$ for $(x, y)\in \R^n\times \R$. Then the optimal value function $\mu$ given in \eqref{mf} has the representation
\begin{equation*}
\mu(x)=\inf_{y\in F(x)}f(x, y)=\inf\{ f_1(x)+y\; |\; y\geq f_2(x)\}=f_1(x)+f_2(x)
\end{equation*} for all $x\in \R^n$.
Note that this representation holds true even if $f_2(x)=\infty$ due to the convention $\inf \emptyset=\infty$. Observe that
\begin{equation*}
\dom f=(\dom f_1)\times \R\ \; \mbox{\rm and }\; \dom F=\epi f_1.
\end{equation*}
Choose $\hat{x}\in \mbox{\rm ri}(\dom f_1)\cap \mbox{\rm ri}(\dom f_2)$ and take $\hat{\lambda}>f_2(\hat{x})$. Then $(\hat{x}, \hat{\lambda})\in \mbox{\rm ri}(\dom f)\cap \mbox{\rm ri}(\gph F)$ by \cite[Proposition~3.7]{NTY2023}. For $w\in \R^n$ and $\gamma\in \R$ we have  $F^*(w, -1)=f^*_2(w)$ and
\begin{equation}\label{fstar}
f^*(w, \gamma)=
\begin{cases}f_1^*(w)&\mbox{\rm if }\gamma=1,\\
 \infty&\mbox{\rm if }\gamma\neq 1.
\end{cases}
\end{equation}
To continue the proof, we observe that
\begin{equation*}
f^*_1(w_1)+f_2^*(w_2)\geq (f_1^*\s f_2^*)(w)\geq (f_1+f_2)^*(w)=\mu^*(w)
\end{equation*}
whenever $w_1+w_2=w$ by the definition. Consider the case where $(f_1+f_2)^*(w)\in \R$. Then $\mu^*(w)=(f_1+f_2)^*(w)\in \R$, so by Theorem~\ref{FC1}, there exist $w_1, w_2\in \R^n$ and $\gamma\in \R$ such that $w=w_1+w_2$ and
\begin{equation*}
\mu^*(w)=f^*(w_1, \gamma)+f^*_2(w_2, -\gamma).
\end{equation*}
Then $\gamma=1$ by \eqref{fstar} and thus
\begin{equation*}
\mu^*(w)=f^*(w_1, 1)+F^*(w_2, -1)=f^*_1(w_1)+f^*_2(w_2)\geq (f_1^*\s f_2^*)(w).
\end{equation*}
This implies equality \eqref{sumf} along with the last statement in this theorem. $\h$

As another direct consequence of Theorem~\ref{FC1}, we obtain a formula for the Fenchel conjugate of the composition of a nearly convex function and an affine mapping.

\begin{Theorem}\label{FC3} Let $g\colon \R^p\to \oR$ be a proper nearly convex function and let $A\in \R^{p\times n}$. Suppose that
\begin{equation*}
A(\R^n)\cap \mbox{\rm ri}(\dom g)\neq\emptyset.
\end{equation*}
Then we have the representation
\begin{equation*}
(g\circ A)^*(w)=\inf\{g^*(v)\; |\; A^Tv=w\}, \ \; w\in \R^n.
\end{equation*}\
In addition, if $(g\circ A)^*(w)\in \R$, then there exists a vector $v\in \R^p$ such that $A^Tv=w$ and $(g\circ A)(w)=g^*(v)$.
\end{Theorem}
{\bf Proof.} Let $f(x, y)=g(y)$  and let $F(x)=\{Ax\}$ for $(x, y)\in \R^n\times \R^p$. Then using the function $\mu$ defined in \eqref{mf}, we have
\begin{equation*}
\mu(x)=(g\circ A)(x), \ \; x\in \R^n,
\end{equation*}
It is not hard to verify that $f^*(u, v)=g^*(v)$ if $u=0$, and $f^*(u, v)=\infty$ otherwise. Moreover,
\begin{equation*}
F^*(u, v)=\begin{cases} 0&\mbox{\rm if }-A^Tv=u,\\
\infty&\mbox{\rm otherwise}.
\end{cases}
\end{equation*}
The conclusion then follows directly from Theorem~\ref{FC1}. $\h$

The example below illustrates the use of Theorem~\ref{FC2} in calculating a useful Fenchel conjugate.
\begin{Example}\label{FC4}{\rm Let $g\colon \R^n\to \oR$ and $h\colon \R^p\to \oR$ be proper nearly convex functions. Given $A\in \R^{p\times n}$, define
\begin{equation*}
f(x, y)=g(x)+h(Ax+y), \ \; (x, y)\in \R^n\times \R^p.
\end{equation*}
Then we have
\begin{equation*}
f^*(0, y^*)=g^*(-A^Ty^*)+h^*(y^*)\; \mbox{\rm for all }y^*\in \R^p.
\end{equation*}
Indeed, let $f_1(x, y)=g(x)$ and $f_2(x, y)=h(Ax+y)$ for $(x, y)\in \R^n\times \R^p$. Fix $(x^*, y^*)\in \R^n\times \R^p$. Then it follows from the definition of the Fenchel conjugate that
\begin{equation*}
f^*_1(x^*, y^*)=
\begin{cases}
g^*(x^*)&\; \mbox{\rm if }y^*=0,\\
\infty &\; \mbox{\rm if }y^*\neq 0.
\end{cases}
\end{equation*}
Observe that $f_2(x, y)=g(\tilde{A}(x, y))$, where $\Tilde{A}(x, y)=Ax+y$. Using Corollary~\ref{FC3}, we have
\begin{equation*}
f_2^*(x^*, y^*)=\inf\{h^*(u^*)\; |\; A^Tu^*=x^*, \; u^*=y^*\}=
\begin{cases} h^*(y^*)\; &\mbox{\rm if }A^Ty^*=x^*,\\
\infty \; &\mbox{\rm if }A^Ty^*\neq x^*.
\end{cases}
\end{equation*}
Next, we use Theorem~\ref{FC2} and get
\begin{equation*}
f^*(0, y^*)=(f_1^*\s f_2^*)(0, y^*)=\inf\{f^*_1(u^*, y^*_1)+f^*_2(-u^*, y^*_2)\; |\; u^*\in \R^n, \; y^*_1+y^*_2=y^*\}.
\end{equation*}
Since $f^*_1(u^*, y^*_1)=\infty$ if $y^*_1\neq 0$, and $f^*_2(-u^*, y^*_2)=\infty$ if $A^Ty^*_2\neq -u^*$, we are only concerned with the case where $y^*_1=0$ and $A^Ty^*_2= -u^*$. In this case, we have $y^*_2=y^*$ and thus $u^*=-A^Ty^*$. It follows that
\begin{equation*}
f^*(0, y^*)=g^*(u^*)+h^*(y^*_2)=g^*(-A^Ty^*)+h^*(y^*)
\end{equation*}
as desired. Note that the qualification conditions in both Theorem~\ref{FC2} and Theorem~\ref{FC3} are satisfied.
}
\end{Example}

\section{Duality Theory for Nearly Convex Optimization}
\setcounter{equation}{0}

In this section, we use the results from the previous sections to develop {\color{blue} a} duality theory for some general optimization models with nearly convex data. Our developments generalize existing results from~\cite{bkw2008}, while providing alternative simpler proofs of these results based on a geometric approach. To develop the duality theory, we follow the approach of using the optimal value function, which has been proven to be successful in the literature; see~\cite{bkw2008,r1,Zalinescu2002}.

Given a function $\ph\colon \R^n\to \oR$, consider the next optimization problem known as the {\em primal problem}:
\begin{equation*}
\mbox{\rm minimize }\; \ph(x), \ \; x\in \R^n.
\end{equation*}
Suppose that $f\colon \R^n\times \R^p\to \oR$ is a  function such that
\begin{equation*}
\ph(x)=f(x, 0)\ \; \mbox{\rm for all } x\in \R^n.
\end{equation*}
Then the problem above can be written as
\begin{equation}\label{P1}
\mbox{\rm minimize }\; f(x, 0), \ \; x\in \R^n.
\end{equation}
Following Rockafellar~\cite{r1}, we consider the following {\em dual problem}:
\begin{equation}\label{D}
\mbox{\rm maximize }-f^*(0, y^*), \ \; y^*\in \R^p.
\end{equation}
Define the optimal values of the primal problem and dual problem, respectively, by
\begin{equation}\label{vvd}
\mathcal{V}=\inf\{f(x, 0)\; |\; x\in \R^n\}\ \; \mbox{\rm and }\; \mathcal{V}_d=\sup\{-f^*(0, y^*)\; |\; y^*\in \R^p\}.
\end{equation}
We also consider the optimal value function
\begin{equation}\label{mff}
\mu(y)=\inf\{f(x, y)\; |\; x\in \R^n\}, \ \; y\in \R^p
\end{equation} of the parametric optimization problem
\begin{equation*}\label{Py}
	\mbox{\rm minimize }\; f(x,y), \ \; x\in \R^n.
\end{equation*}
For a function $h\colon \R^n\to \oR$, the second Fenchel conjugate of $h$ is defined by $h^{**}=(h^*)^*$.

The theorem below shows that $\mathcal{V}$ and $\mathcal{V}_d$ can be expressed in term of the optimal value function~\eqref{mff} and its Fenchel second conjugate.
\begin{Theorem}\label{SD1} {\rm  (cf.~\cite[Theorem~2.6.1]{Zalinescu2002})} Consider the function $\mu$ given in~\eqref{mff} in which $f\colon \R^n\times \R^p\to \oR$ is a  nearly convex function. Then we have the following assertions:
\begin{enumerate}
\item $\mathcal{V}=\mu(0)$ and $\mathcal{V}_d=\mu^{**}(0)$.
\item If $\partial \mu(0)\neq\emptyset$ with $\mu(0)\in \R$, then $\mathcal{V}=\mathcal{V}_d$.
\item If $\mu$ is proper, convex, and lower semicontinuous, then $\mathcal{V}=\mathcal{V}_d$.
\end{enumerate}
\end{Theorem}
{\bf Proof.} It follows from the definition that
\begin{equation*}
\mu(0)=\inf\{f(x, 0)\; |\; x\in \R^n\}=\mathcal{V}.
\end{equation*}
In addition,
\begin{equation*}
\mu^{**}(0)=\sup\{\la 0, y^*\ra-\mu^*(y^*)\; |\; y^*\in \R^n\}=\sup\{-\mu^*(y^*)\; |\; y^*\in \R^n\}.
\end{equation*}
Thus, by Corollary~\ref{FC1.1} and Remark~\ref{FC1.12} we have
\begin{equation*}
\mu^{**}(0)=\sup\{-\mu^*(y^*)\; |\; y^*\in \R^n\}=\sup\{-f^*(0, y^*)\; |\; y^*\in \R^n\}=\mathcal{V}_d.
\end{equation*}
This completes the proof of (a).

Now, suppose that $\partial \mu(0)\neq\emptyset$. By Young's inequality, for $y, z\in \R^p$ we have
\begin{equation*}
\la y, z^*\ra\leq \mu(y)+ \mu^*(z^*),
\end{equation*}
which implies
\begin{equation*}
\la y, z^*\ra-\mu^*(z^*)\leq \mu(y).
\end{equation*}
Taking the supremum from both sides of this inequality yields $\mu^{**}(y)\leq \mu(y)$, and in particular, $\mu^{**}(0)\leq \mu(0)$. Pick $z\in \partial \mu(0)$. Then
\begin{equation*}
\la z^*, y\ra - \mu(y)\leq -\mu(0)\ \; \mbox{\rm for all }y\in \R^p.
\end{equation*}
Thus, taking the supremum with respect to $y\in \R^p$ gives us $\mu^*(z^*)\leq -\mu(0)$. Now, we have
\begin{equation*}
\mu^{**}(0)=\sup\{\la 0, y\ra-\mu^*(y)\; |\; y\in \R^n\}=\sup\{-\mu^*(y^*)\; |\; y^*\in \R^n\}\geq -\mu^*(z^*)\geq \mu(0),
\end{equation*}
which justifies the equality $\mu(0)=\mu^{**}(0)$ and also assertion (b). If $\mu$ is proper, convex, and lower semicontinuous, then $\mu=\mu^{**}$ by \cite[Theorem~4.10]{bmn}. Therefore, assertion (c) follows from (a). $\h$

As a corollary, we obtain the following result in \cite{bkw2008}.
\begin{Corollary}\label{GD} Let $f\colon \R^n\times \R^p\to \oR$ be a proper nearly convex function. Consider the primal problem \eqref{P1} and its dual problem \eqref{D} along with their optimal values $\mathcal{V}$ and $\mathcal{V}_d$ given in \eqref{vvd}.  Suppose $0\in \mbox{\rm ri}(\mathcal{P}_2(\dom f))$, where $\mathcal{P}_2\colon \R^n\times \R^p\to \R^p$ is the projection mapping defined by
	\begin{equation}\label{P2}
		\mathcal{P}_2(x, y)=y, \ \; (x, y)\in \R^n\times \R^p.
	\end{equation} Then $\mathcal{V}=\mathcal{V}_d$.
\end{Corollary}
{\bf Proof.} Consider the  function $\mu$ given in \eqref{mff}.  Obviously, $\dom \mu=\mathcal{P}_2(\dom f)$.
Observe that $\mu^{**}(0)\leq \mu(0)$ as in the proof of Theorem~\ref{SD1}. If $\mu(0)=-\infty$, then $\mu^{**}(0)=-\infty$ and thus $\mathcal{V}$ and $\mathcal{V}_d$ in this case by   Theorem~\ref{SD1}(a). Consider the case where $\mu(0)>-\infty$. Since $0\in \mbox{\rm ri}(\mathcal{P}_2(\dom f))=\mbox{\rm ri}(\dom \mu)$, we see by  \cite[Proposition~5.4]{NTY2023} that $\partial \mu(0)\neq\emptyset$. Note that $\mu$ is a proper function by Proposition~\ref{PP}. This implies by Theorem \ref{SD1} that $\mathcal{V}=\mu(0)=\mu^{**}(0)=\mathcal{V}_d$. $\h$

Now, we consider the problem
\begin{eqnarray}\label{primal p}
\begin{array}{ll}
&\mbox{\rm minimize } \phi(x)\\
&\mbox{\rm subject to }\; 0\in G(x), \; x\in \Theta,
\end{array}
\end{eqnarray}
where $\phi\colon \R^n\to \oR$,  $\Theta\subset\R^n$, and $G\colon \R^n\tto \R^q$ is a set-valued mapping.  This problem can be converted to minimizing the function
\begin{equation}\label{phf}
\ph(x)=\phi(x)+\delta(x; G^{-1}(0)\cap \Theta), \ \; x\in \R^n.
\end{equation}
Consider the function $f\colon \R^n\times \R^q\to \oR$ given by
\begin{equation}\label{pt11}
f(x, y)=\begin{cases}\phi(x)\;& \mbox{\rm if }y\in G(x)\; \mbox{\rm and }x\in \Theta,\\
\infty\; &\mbox{\rm otherwise}.
\end{cases}
\end{equation}
It is obvious that $\ph(x)=f(x, 0)$. Consider the optimal value function
\begin{equation}\label{vg}
v_G(x, y^*)=\inf \{\la -y^*, y\ra\; |\; y\in G(x)\},\ \, x\in \R^n\; \mbox{\rm and } y^*\in \R^q.
\end{equation}
The dual problem of  \eqref{primal p} is given by
\begin{eqnarray}\label{dual p}
	\begin{array}{ll}
&\mbox{\rm maximize }h(y^*)=\inf \{\phi(x)+v_G(x, y^*)\; |\; x\in \Theta\},\\
&\mbox{\rm subject to }y^*\in \R^q.
\end{array}
\end{eqnarray}
The theorem below establishes a strong duality result for problem~\eqref{primal p} and its dual problem~\eqref{dual p}.
\begin{Theorem}\label{LD11} Suppose that $\phi\colon \R^n\to \oR$ is a proper nearly convex function, $\Theta$ is a nearly convex set in $\R^n$, and $G\colon \R^n\tto \R^q$ is a nearly convex set-valued mapping.  Consider the problem primal problem \eqref{primal p} and its dual problem~\eqref{dual p}. Suppose that
\begin{equation}\label{QCD1}
\mbox{\rm ri}(\dom \phi)\cap \mbox{\rm ri}(\dom G)\cap \ri \Theta\neq\emptyset\;\ \mbox{\rm and }\; 0\in \mbox{\rm ri}\big(G(\Theta\cap \dom \phi)\big).
\end{equation}
Then we have the Lagrange strong duality, i.e., $\mathcal{V}=\mathcal{V}_d$, where
\begin{align*}
&\mathcal{V}=\inf\{\phi(x)\; |\; x\in \Theta, \; 0\in G(x)\},\\
&\mathcal{V}_d=\sup_{y^*\in \R^q}h(y^*).
\end{align*}
\end{Theorem}
{\bf Proof.} Observe by Corollary~\ref{bot3} that the function $f$ defined in \eqref{pt11} is nearly convex under the first condition in \eqref{QCD1}. Using Theorem~\ref{ri_of_image}, we can also verify that $0\in \mbox{\rm ri}(\mathcal{P}_2(\dom f))$, where $\mathcal{P}_2$ is defined in \eqref{P2}, under the second condition in \eqref{QCD1}.  We have the following equalities:
\begin{align*}
-f^*(0, y^*)&=-\sup\{\la (0, y^*), (x, y)\ra-f(x, y)\; |\; (x, y)\in \R^n\times \R^q\}\\
&=-\sup\{\la y^*, y\ra-\phi(x)\; |\; x\in \Theta, \; y\in G(x)\}\\
&=\inf\{\phi(x)-\la y^*, y\ra\; |\; x\in \Theta, \; y\in G(x)\}\\
&=\inf \{\phi(x)+ v_G(x, y^*)\; |\; x\in \Theta\}=h(y^*).
\end{align*}
Thus, the conclusion is a direct consequence of Corollary~\ref{GD}.  $\h$

\begin{Lemma}\label{inc} Suppose that $G(x)=g(x)+K$, $x\in \R^n$, where $g\colon\R^n\to \R^q$ and $K$ is a nonempty convex cone. Then the function $v_G$ given in \eqref{vg} has the following representation:
\begin{equation}\label{vrep}
v_G(x, y^*)=\begin{cases} -\la y^*, g(x)\ra\; &\mbox{\rm if }-y^*\in K^*,\\
-\infty\; &\mbox{\rm otherwise}.
\end{cases}
\end{equation}
where $K^*=\{y^*\in \R^q\; |\; \la y^*, y\ra\geq 0\; \mbox{\rm for all }\, y\in K\}$.
\end{Lemma}
{\bf Proof.} Observe that
\begin{align*}
v_G(x, y^*)=\inf\{\la -y^*, y\ra\; |\; y\in G(x)\}&=\inf\{\la -y^*, y\ra\; |\; y\in g(x)+K\}\\
&=\la -y^*, g(x)\ra+\inf\{\la -y^*, y\ra\; |\; y\in K\}.
\end{align*}
If $-y^*\in K^*$, then $\la -y^*, y\ra\geq 0$ for all $y\in K$, and $\la -y^*, 0\ra=0$ with $0\in K$, so $\inf\{\la -y^*, y\ra\; |\; y\in K\}=0$. It follows that
$v_G(x, y^*)=-\la y^*, g(x)\ra$. In the case where $-y^*\notin K^*$, we can find $y_0\in K$ such that $\la -y^*, y_0\ra<0$, so $$\inf\{\la -y^*, y\ra\; |\; y\in K\}\leq \inf\{\la -y^*, ty_0\ra\; |\; t>0\}=-\infty$$ because $K$ is a cone. Thus, we obtain \eqref{vrep}. $\h$

Next, consider the problem
\begin{eqnarray}\label{BotP}
\begin{array}{ll}
&\mbox{\rm minimize } \phi(x)\\
&\mbox{\rm subject to }g(x)\leq_K 0, \; x\in \Theta.
\end{array}
\end{eqnarray}
We will study the relationship between problem \eqref{BotP} and its dual problem:
\begin{eqnarray}\label{BotD}
\begin{array}{ll}
&\mbox{\rm maximize } \ph(y^*)=\inf\{\phi(x)+\la y^*, g(x)\ra\; |\; x\in \Theta\}\\
&\mbox{\rm subject to }y^*\in K^*.
\end{array}
\end{eqnarray}
Now, we are ready to obtain a result in \cite{bkw2008}.
\begin{Corollary}\label{B duality} {\rm\bf (\cite[Theorem~4.2]{bkw2008})} Suppose that $\phi\colon \R^n\to \oR$ is a proper nearly convex function,  $\Theta$ is a nearly convex set in $\R^n$, and  $g\colon \R^n\to \R^q$ is a $K$-nearly convex, where $K$ is a nonempty convex cone. Consider the primal problem \eqref{BotP} along with its dual problem~\eqref{BotD}. Suppose that
\begin{equation*}
\ri\Theta\cap \mbox{\rm ri}(\dom \phi)\neq\emptyset\ \; \mbox{\rm and }\; 0\in g(\mbox{\rm ri}(\Theta)\cap \mbox{\rm ri}(\dom \phi))+\ri K.
\end{equation*}
Then we have the Lagrange strong duality, i.e., $\mathcal{V}=\mathcal{V}_d$, where
\begin{align*}
&\mathcal{V}=\inf\{\phi(x)\; |\; g(x)\leq_K 0, \; x\in \Theta\},\\
&\mathcal{V}_d=\sup\{\ph(y^*)\; |\; y^*\in K^*\}.
\end{align*}
\end{Corollary}
{\bf Proof.} Consider the set-valued mapping $G(x)=g(x)+K$ for $x\in \R^n$ and get $\dom G=\R^n$. Using Lemma~\ref{inc}, we see that the function $h$ in Theorem~\ref{LD11} has the representation
\begin{equation*}
h(y^*)=\begin{cases}
\inf\{\phi(x)+\la -y^*, x\ra\; |\; x\in \Theta\}\; &\mbox{\rm if }-y^*\in K^*,\\
-\infty \; &\mbox{\rm otherwise}.
\end{cases}
\end{equation*}
It follows that
\begin{align*}
\sup_{y^*\in \R^q}h(y^*)=\sup_{y^*\in -K^*}\inf_{x\in \Theta}\, [\phi(x)+\la -y^*, x\ra] &=\sup_{y^*\in K^*}\inf_{x\in \Theta}\, [\phi(x)+\la y^*, x\ra]\\ & =\sup\{\ph(y^*)\; |\; y^*\in K^*\}.
\end{align*}}
Finally, we check that the qualification conditions in Theorem~\ref{LD11} are satisfied. Since $\dom G=\R^n$, it is obvious that
\begin{equation*}
\mbox{\rm ri}(\dom G)\cap \ri \Theta\cap \mbox{\rm ri}(\dom \phi)\neq\emptyset.
\end{equation*}
The condition $0\in g(\mbox{\rm ri}(\Theta)\cap \mbox{\rm ri}(\dom \phi))+\ri K$ guarantees that
\begin{equation*}
0\in \mbox{\rm ri}\big(G(\Theta\cap \dom \phi)\big)
\end{equation*}
by using Corollary~\ref{ri epi f}. Thus, the conclusion follows  from Theorem~\ref{LD11}. $\h$

Next, we consider problem \eqref{primal p} again with the perturbed function $f_1\colon \R^n\times \R^n\times \R^q\to \oR$ given by
\begin{equation}\label{pt2}
f_1(x, u, y)=\begin{cases}\phi(x+u)\;& \mbox{\rm if }u\in \R^n, \; y\in G(x)\; \mbox{\rm and }x\in \Theta,\\
\infty\; &\mbox{\rm otherwise}.
\end{cases}
\end{equation}
It is obvious that $\ph(x)=f_1(x, 0, 0)$, where $\ph$ is given in \eqref{phf}.

We will study  the  primal problem \eqref{primal p} in connection with its dual problem:
\begin{eqnarray}\label{dual p1}
\begin{array}{ll}
&\mbox{\rm maximize }h_1(u^*, y^*)=-\phi^*(u^*)+\inf \{\la u^*, x\ra +v_G(x, y^*)\; |\; x\in \Theta\},\\
&\mbox{\rm subject to }u^*\in \R^n, \; y^*\in \R^q.
\end{array}
\end{eqnarray}
Based on this perturbation, the next theorem gives us a new strong duality result.

\begin{Theorem}\label{LD1} Suppose that $\phi\colon \R^n\to \oR$ is a proper nearly convex function, $\Theta$ is a nearly convex set in $\R^n$, and  $G\colon \R^n\tto \R^q$ is a nearly convex set-valued mapping.  Consider the problem primal problem \eqref{primal p} and it dual problem \eqref{dual p1}. Suppose that \eqref{QCD1} is satisfied. Then we have the Lagrange-Fenchel strong duality, i.e.,  $\mathcal{V}=\mathcal{V}_d$, where
\begin{align*}
&\mathcal{V}=\inf\{\phi(x)\; |\; x\in \Theta, \; 0\in G(x)\},\\
&\mathcal{V}_d=\sup_{u^*\in \R^n, \; y^*\in \R^q}h_1(u^*, y^*).
\end{align*}
\end{Theorem}
{\bf Proof.} By Corollary~\ref{bot33}, the function $f_1$ given in~\eqref{pt2} is nearly convex under the first condition in \eqref{QCD1}. Fix any $u^*\in \R^n$ and $y^*\in \R^q$. Let us  show that
\begin{equation*}
-f_1^*(0, u^*, y^*)=h_1(u^*, y^*)\; \mbox{\rm for all }(u^*, y^*)\in \R^n\times \R^q.
\end{equation*}
Indeed, using $z=x+u$, we have
\begin{align*}
f_1^*(0, u^*, y^*)=&\sup\{\la u^*, u\ra +\la y^*, y\ra-f_1(x, u, y)\; |\; (x, u, y)\in \R^n\times \R^n\times \R^q\}\\
&=\sup\{\la u^*, u\ra+\la y^*, y\ra -\phi(x+u)\; |\; x\in \Theta, \; y\in G(x), \; u\in \R^n\}\\
&=\sup\{\la u^*, z-x\ra+\la y^*, y\ra-\phi(z)\; |\; x\in \Theta, \; y\in G(x), \; z\in \R^n\}\\
&=\sup\{\la u^*, z\ra-\phi(z)\; |\; z\in \R^n\}+\sup\{-\la u^*, x\ra +\la y^*, y\ra\; |\; x\in \Theta, \; y\in G(x)\}\\
&=\phi^*(u^*)-\inf\{\la u^*, x\ra+\la -y^*, y\ra\; |\; x\in \Theta, \; y\in G(x)\}\\
&=\phi^*(u^*)-\inf\{\la u^*, x\ra +v_G(x, y^*)\; |\; x\in \Theta\}.
\end{align*}
It follows that
\begin{equation*}
-f^*_1(0, u^*, y^*)=-\phi^*(u^*)+\inf\{\la u^*, x\ra +v_G(x, y^*)\; |\; x\in \Theta\}=h_1(u^*, y^*).
\end{equation*}
The second condition in \eqref{QCD1} with the use of Theorem~\ref{ri_of_image} guarantees that there exists $\hat{x}\in \mbox{\rm ri}(\dom G)\cap \ri\Theta\cap \mbox{\rm ri}(\dom \phi)$ such that $0\in \ri G(\hat{x})$. Then the element $(\hat{x}, \hat{u}, \hat{y})$, where $\hat{u}=0$ and $\hat{y}=0$, belongs to the set $\mbox{\rm ri}(\dom  f_1)$. It follows that $(0, 0)\in \mbox{\rm ri}\big(\widetilde{\mathcal{P}}_{2}(\dom f_1)\big)$, where $\widetilde{\mathcal{P}}_2\colon \R^n\times \R^n\times \R^p$ is given by
\begin{equation*}
\widetilde{\mathcal{P}}_2(x, u, y)=(u, y), \ \; (x, u, y)\in \R^n\times \R^n\times \R^q.
\end{equation*}
Therefore, the conclusion follows from Corollary~\ref{GD}. $\h$

\begin{Remark}{\rm  Using Theorem~\ref{LD1} with $G(x)=g(x)+K$ as in the setting of Corollary~\ref{B duality} allows ones to obtain \cite[Theorem~4.4]{bkw2008}. We leave this as an exercise for the reader. }
\end{Remark}

Let $g\colon \R^n\to \oR$ and $h\colon \R^p\to \oR$ be proper nearly convex functions. Given $A\in \R^{p\times n}$, consider the problem
\begin{equation}\label{FP}
\mbox{\rm minimize } \ph(x)=g(x)+h(Ax), \ \; x\in \R^n.
\end{equation}
Then $\ph(x)=f(x, 0)$, where
\begin{equation}\label{Fenchelf}
f(x, y)=g(x)+h(Ax+y), \ \; (x, y)\in \R^n\times \R^p.
\end{equation}
By Example~\ref{FC4}, the dual problem \eqref{D} becomes
\begin{equation}\label{FD}
\mbox{\rm maximize }-f^*(0, y^*)=-g^*(-A^Ty^*)-h^*(y^*), \ \; y^*\in \R^p.
\end{equation}
We are ready to recover a result in \cite{bkw2008} on Fenchel duality.

\begin{Theorem}{\rm\bf (\cite[Theorem~4.6]{bkw2008})} Let $g\colon \R^n\to \oR$, $h\colon \R^p\to \oR$ be proper nearly convex functions, and let $A\in \R^{p\times n}$. Consider the primal problem~\eqref{FP} and the dual problem~\eqref{FD}. Suppose that
\begin{equation*}
A\big(\mbox{\rm ri}(\dom g)\big)\cap \mbox{\rm ri}(\dom h)\neq\emptyset.
\end{equation*}
Then the Fenchel strong duality holds, i.e., $\mathcal{V}=\mathcal{V}_d$, where
\begin{align*}
&\mathcal{V}=\inf \{g(x)+h(Ax)\; |\; x\in \R^n\},\\
&\mathcal{V}_d=\sup\{-g^*(-A^Ty^*)-h^*(y^*)\; |\; y^*\in \R^p\}.
\end{align*}
\end{Theorem}
{\bf Proof.} By Corollary~\ref{GD}, it suffices to show that $0\in \mbox{\rm ri}(\mathcal{P}_2(\dom f))$, where $f$ is defined in~\eqref{Fenchelf} and $\mathcal{P}_2$ is given by~\eqref{P2}. Choose $\hat{x}\in \mbox{\rm ri}(\dom g)$ such that $A\hat{x}\in \mbox{\rm ri}(\dom h)$. Then set $\hat{y}=0$ and get $f(\hat{x}, \hat{y})=g(\hat x)+h(A\hat x)\in \R$, so $(\hat{x}, 0)\in \mbox{\rm ri}(\dom f)$. This implies that $0\in \mbox{\rm ri}(\mathcal{P}_2(\dom f))$ and completes the proof. $\h$
\\[1ex]
{\bf Acknowledgment.} Nguyen Mau Nam would like to thank  the Vietnam Institute for Advanced Study in Mathematics for hospitality.

\end{document}